\documentclass[aos,MSNbibl,nameyear,secthm,seceqn,dvips]{arximspdf}
\usepackage{dcolumn}
\usepackage{graphicx}

\doi{10.1214/10-AOS871}
\volume{39}
\issue{3}
\pubyear{2011}
\firstpage{1658}
\lastpage{1688}

\makeatletter
\newcolumntype{d}[1]{D{.}{.}{#1}}

\newproclaim{rem}{Remark}
\newproclaim{exm}{Example}
\newtheorem{cor}{Corollary}
\newtheorem{prop}{Proposition}

\newcommand{\CO}{{\mathcal{O}}}
\newcommand{\R}{\mathbb{R}}

\renewcommand{\citep}[1]{[\citet{#1}]}

\newcommand{\mo}{\mbox{\fontsize{8}{9}\selectfont{$\mathcal{O}$}}}

\let\sv@tabnotetext\tabnotetext
\let\sv@tabnotemark@fmt\tabnotemark@fmt
\long\def\legend#1{{\let\tabnote@indent\leavevmode\sv@tabnotetext[]{}{#1}}}

\setattribute{abstract}{width}{290pt}

\makeatother

\begin{document}
\begin{frontmatter}

\title{The EFM approach for single-index models}
\runtitle{The EFM approach for single-index models}

\begin{aug}
\author[A]{\fnms{Xia} \snm{Cui}\thanksref{t1}\ead[label=e1]{cuixia@mail.sysu.edu.cn}},
\author[B]{\fnms{Wolfgang Karl} \snm{H\"{a}rdle}\thanksref{t2}\ead[label=e2]{haerdle@wiwi.hu-berlin.de}}
and
\author[C]{\fnms{Lixing} \snm{Zhu}\corref{}\thanksref{t3}\ead[label=e3]{lzhu@hkbu.edu.hk}}
\runauthor{X. Cui, W. K. H\"{a}rdle and L. Zhu}
\affiliation{Sun Yat-sen University,
Humboldt-Universit\"{a}t zu Berlin and\\ National Central University,
and Hong
Kong Baptist University\\ and Yunnan University of Finance and Economics}
\thankstext{t1}{Supported by NNSF project (11026194)
of China, RFDP (20100171120042) of China and ``the Fundamental Research
Funds for the Central Universities'' (11lgpy26) of China.}
\thankstext{t2}{Supported by Deutsche
Forschungsgemeinschaft SFB 649 ``\"{O}konomisches Risiko.''}
\thankstext{t3}{Supported by a Grant (HKBU2030/07P) from Research Grants
Council of Hong Kong, Hong Kong, China.}
\address[A]{X. Cui\\
School of Mathematics \\
\quad and Computational Science\\
Sun Yat-sen University\\
Guangzhou\\
Guangdong Province, 510275\\
P.R. China\\
\printead{e1}}

\address[B]{W. K. H\"{a}rdle\\
CASE-Center for Applied Statistics\\
\quad and Economics\\
Humboldt-Universit\"{a}t zu Berlin\\
Wirtschaftswissenschaftliche Fakult\"{a}t\\
Spandauer Str. 1\\
10178 Berlin \\
Germany\\
\printead{e2}
}

\address[C]{L. Zhu\\
FSC1207, Fong Shu Chuen Building\\
Department of Mathematics\\
Hong Kong Baptist University\\
Kowloon Tong\\
Hong Kong\\
P.R. China\\
\printead{e3}}
\end{aug}

\received{\smonth{4} \syear{2010}}
\revised{\smonth{12} \syear{2010}}

%
\begin{abstract}
Single-index models are natural extensions of linear models and
circumvent the so-called curse of dimensionality. They are becoming
increasingly popular in many scientific fields including biostatistics,
medi\-cine, economics and financial econometrics. Estimating and testing
the model index coefficients $\bolds{\beta}$ is one of the most important
objectives in the statistical analysis. However, the commonly used
assumption on the index coefficients, $\|\bolds{\beta}\|=1$,
represents a
nonregular problem: the true index is on the boundary of the unit
ball. In this paper we introduce the EFM approach, a method of
estimating functions, to study the single-index model. The procedure is
to first relax the equality constraint to one with $(d-1)$ components
of $\bolds{\beta}$ lying in an open unit ball, and then to construct the
associated $(d-1)$ estimating functions by projecting the score
function to the linear space spanned by the residuals with the unknown
link being estimated by kernel estimating functions. The root-$n$
consistency and asymptotic normality for the estimator obtained from
solving the resulting estimating equations are achieved, and a Wilks
type theorem for testing the index is demonstrated. A noticeable result
we obtain is that our estimator for $\bolds{\beta}$ has smaller or equal
limiting variance than the estimator of Carroll et al. [\textit{J. Amer. Statist. Assoc.}
\textbf{92} (1997) 447--489]. A fixed-point iterative scheme for
computing this estimator is proposed.
This algorithm only involves one-dimensional nonparametric smoothers,
thereby avoiding the data sparsity problem caused by high model
dimensionality. Numerical studies based on simulation and on
applications suggest that this new estimating system is quite powerful
and easy to implement.
\end{abstract}

%


\end{frontmatter}

\section{Introduction}\label{sec1}

Single-index models combine flexibility of modeling with
interpretability of (linear) coefficients.
They circumvent the curse of dimensionality and
are becoming increasingly popular in many scientific fields. The
reduction of
dimension is achieved by assuming the link function to be a univariate
function applied to the projection of explanatory covariate vector on
to some direction. In this paper we consider an extension of
single-index models where, instead of a distributional assumption,
assumptions of only the mean function and variance function of the
response are made. Let $(Y_i,{\mathbf{X}}_i)$, $i=1,\ldots,n,$ denote
the observed values with $Y_i$ being the response variable and ${\mathbf{X}}_i$ as
the vector of $d$ explanatory variables. The relationship of the mean
and variance of ${ Y}_i$ is specified as follows:
%
\begin{equation}\label{model}
E(Y_i|{\mathbf{X}}_i)=\mu\{g({\bolds{\beta}}^{\top}{\mathbf{X}}_i)\},\qquad
\operatorname{Var}(Y_i|{\mathbf{X}}_i)=\sigma^2V\{g({\bolds{\beta}}^{\top}{\mathbf{X}}_i)\},
\end{equation}
where $\mu$ is a known monotonic function, $V$ is a known covariance
function, $g$ is an unknown univariate link function and $\bolds{\beta}$
is an unknown index vector which belongs to the parameter space
$\Theta=\{{\bolds{\beta}}=(\beta_1,\ldots,\beta_d)^{\top}\dvtx
\|{\bolds{\beta}}\|=1, \beta_1>0, {\bolds{\beta}}\in{\R}^d\}$. Here we assume
the parameter space is $\Theta$ rather than the entire ${\R}^d$ in
order to ensure that $\bolds{\beta}$ in the representation (\ref{model})
can be uniquely defined. This is a commonly used assumption on the
index parameter [see \citet{Carroll1997}, \citet{Zhu2006}, \citet{Linwei2007}].
Another reparameterization is to let $\beta_1=1$ for
the sign identifiability and to transform ${\bolds{\beta}}$ to
$(1,\beta_2,\ldots,\beta_d)/(1+\sum_{r=2}^d\beta_r^2)^{1/2}$
for the scale identifiability. Clearly
$(1,\beta_2,\ldots,\beta_d)/(1+\sum_{r=2}^d\beta_r^2)^{1/2}$
can also span the parameter space $\Theta$ by simply checking that
$\|(1,\beta_2,\ldots,\beta_d)/(1+\sum_{r=2}^d\beta_r^2)^{1/2}\|=1$
and the first component
$1/(1+\sum_{r=2}^d\beta_r^2)^{1/2}>0$. However, the
fixed-point algorithm recommended in this paper for normalized
vectors may not be suitable for such a reparameterization. Model
(\ref{model}) is flexible enough to cover a variety of situations.
If $\mu$ is the identity function and $V$ is equal to constant~1,
(\ref{model}) reduces to a single-index model \citet{Hardle1993}.
Model (\ref{model}) is an extension of the generalized linear model
\citet{McCullagh1989} and the single-index model. When the\vspace*{1pt}
conditional distribution of $Y$ is logistic, then
$\mu\{g({\bolds{\beta}}^{\top}{\mathbf{X}})\}= \exp\{g({\bolds{\beta}}^{\top
}{\mathbf{X}})\}/[1+\exp\{g({\bolds{\beta}}^{\top}{\mathbf{X}})\}]$ and\vspace*{1pt}
$V\{g({\bolds{\beta}}^{\top}{\mathbf{X}})\}=\exp\{g({\bolds{\beta}}^{\top}{\mathbf{X}})\}/[1+\exp\{g({\bolds{\beta}}^{\top}{\mathbf{X}})\}]^2$.

For single-index models: $\mu\{g({\bolds{\beta}}^{\top}{\mathbf{X}})\}=g({\bolds{\beta}}^{\top}{\mathbf{X}})$ and $ V\{g({\bolds{\beta}}^{\top
}{\mathbf{X}})\}=1$, various strategies for estimating $\bolds{\beta}$ have been
proposed in the last decades. Two most popular methods are the
average derivative method (ADE) introduced in \citet{Powell1989}
and \citet{Hardle1989}, and the simultaneous minimization method
of \citet{Hardle1993}. Next we will review these two methods in
short. The ADE method is based on that $\partial E(Y|{\mathbf{X}}={\mathbf{x}})/\partial{\mathbf{x}}=g'({\bolds{\beta}}^{\top}{\mathbf{x}}){\bolds{\beta}}$ which
implies that the gradient of the regression function is proportional
to the index parameter $\bolds{\beta}$. Then a natural estimator for
$\bolds{\beta}$ is $\hat{\bolds{\beta}}=n^{-1}\sum_{i=1}^n
\widehat{\nabla G}({\mathbf{X}}_i)/\|n^{-1}\sum_{i=1}^n
\widehat{\nabla G}({\mathbf{X}}_i)\|$ with $\nabla G({\mathbf{x}})$ denoting
$\partial E(Y|{\mathbf{X}}={\mathbf{x}})/\partial{\mathbf{x}}$ and $\|\cdot\|$
being the Euclidean norm. An advantage of the ADE approach is that
it allows estimating $\bolds{\beta}$ directly. However, the
high-dimensional kernel smoothing used for computing
$\widehat{\nabla G}({\mathbf{x}})$ suffers from the ``curse of
dimensionality'' if the model dimension $d$ is large.
\citet{Hristache2001} improved the ADE approach by lowering the
dimension of the kernel gradually. The method of
\citet{Hardle1993} is carried out by minimizing a least squares
criterion based on nonparametric estimation of the link $g$ with
respect to $\bolds{\beta}$ and bandwidth~$h$. However, the minimization
is difficult to implement since it depends on an optimization
problem in a high-dimensional space. \citet{Xia2002} proposed to
minimize average conditional variance (MAVE). Because the kernel
used for computing $\bolds{\beta}$ is a function of $\|{\mathbf{X}}_i-{\mathbf{X}}_j\|$, MAVE meets the problem of data sparseness. All the above
estimators are consistent under some regular conditions. Asymptotic
efficiency comparisons of the above methods have been discussed in
\citet{Xia2006} resulting in the MAVE estimator of $\bolds{\beta}$
having the same limiting variance as the estimators of
\citet{Hardle1993}, and claiming alternative versions of the
ADE method having larger variance. In addition, \citet{Yuyan2002}
fitted the partially linear single-index models using a penalized
spline method. \citet{HuhPark2002} used the local polynomial
method to fit the unknown function in single-index models. Other
dimension reduction methods that were recently developed in the
literature are sliced inverse regression, partial least squares and
canonical correlation method. These methods handle high-dimensional
predictors; see \citeauthor{ZhuZhua2009} (\citeyear{ZhuZhua2009}, \citeyear{ZhuZhub2009}) and
\citet{ZhouHe2008}.

The main challenges of estimation in the semiparametric model
(\ref{model}) are that the support of the infinite-dimensional
nuisance parameter $g(\cdot)$ depends on the finite-dimensional
parameter $\bolds{\beta}$, and the parameter $\bolds{\beta}$ is on the
boundary of a unit ball. For estimating $\bolds{\beta}$ the former
challenge forces us to deal with the infinite-dimensional nuisance
parameter $g$. The latter one represents a nonregular problem. The
classic assumptions about asymptotic properties of the estimates for
$\bolds{\beta}$ are not valid. In addition, as a model proposed for
dimension reduction, the dimension $d$ may be very high and one
often meets the problem of computation. To attack the above
problems, in this paper we will develop an estimating function
method (EFM) and then introduce a computational algorithm to solve
the equations based on a fixed-point iterative scheme. We first
choose an identifiable parameterization which transforms the boundary
of a unit ball in ${\R}^d$ to the interior of a unit ball in
${\R}^{d-1}$. By eliminating $\beta_1$, the parameter space\vspace*{1pt} $\Theta$
can be rearranged to a form $ \{((1-\sum_{r=2}^d
\beta_r^2)^{1/2},\beta_2,\ldots,\beta_d)^{\top}\dvtx \sum_{r=2}^d
\beta_r^2<1\}$. Then the derivatives\vspace*{1pt} of a function with respect to
$(\beta_2,\ldots,\beta_d)^{\top}$ are readily obtained by the chain
rule and the classical assumptions on the asymptotic normality hold
after transformation. The estimating functions (equations) for
$\bolds{\beta}$ can be constructed by replacing $g({\bolds{\beta}}^{\top
}{\mathbf{X}})$ with $\hat{g}({\bolds{\beta}}^{\top}{\mathbf{X}})$. The estimate
$\hat{g}$ for the nuisance parameter $g$ is obtained using kernel
estimating functions and the smoothing parameter $h$ is selected
using $K$-fold cross-validation. For the problem of testing the
index, we establish a quasi-likelihood ratio based on the proposed
estimating functions and show that the test statistics
asymptotically follow a $\chi^2$-distribution whose degree of
freedom does not depend on nuisance parameters, under the null
hypothesis. Then a Wilks type theorem for testing the index is
demonstrated.

The proposed EFM technique is essentially a unified method of handling
different types of data situations including categorical response
variable and
discrete explanatory covariate vector. The main results of this
research are as
follows:
\begin{enumerate}[(a)]
\item[(a)]\textit{Efficiency}. A surprising result we obtain is that our EFM estimator
for $\bolds{\beta}$ has smaller or equal limiting variance than the
estimator of \citet{Carroll1997}.
\item[(b)]\textit{Computation}. The estimating function system only
involves one-dimen\-sional nonparametric smoothers, thereby avoiding
the data sparsity problem caused by high model dimensionality.
Unlike the quasi-likelihood inference \citep{Carroll1997} where
the maximization is difficult to implement when $d$ is large, the
reparameterization and the explicit formulation of the estimating
functions facilitate an efficient computation algorithm. Here we
use a fixed-point iterative scheme to compute the resultant
estimator. The simulation results show that the algorithm adapts
to higher model dimension and richer data situations than the MAVE
method of \citet{Xia2002}.

\end{enumerate}

It is noteworthy that the EFM approach proposed in this paper cannot
be obtained from the SLS method proposed in \citet{Ichimura1993} and
investigated in
\citet{Hardle1993}. SLS minimizes the weighted least squares
criterion $\sum_{j=1}^n
[Y_j-\mu\{\hat{g}({\bolds{\beta}}^{\top}{\mathbf{X}}_j)\}]^2V^{-1}\{\hat{g}({\bolds{\beta}}^{\top}{\mathbf{X}}_j)\}$, which
leads to a biased estimating equation when we use its derivative if
$V(\cdot)$ does not contain the parameter of interest. It will not
in general provide a consistent estimator [see \citet{Heyde1997}, page 4]. \citet{Chang2010} and
\citet{WangXueZhuChong2010} discussed the efficient estimation
of single-index model for the case of additive noise. However, their
methods are based on the estimating equations induced from the
least squares rather than the quasi-likelihood. Thus, their
estimation does not have optimal property. Also their comparison is
with the one from \citet{Hardle1993} and its later development. It
cannot be applied to the setting under study. In this paper, we
investigate the efficiency and computation of the estimates for the
single-index models, and systematically develop and
prove the asymptotic properties of EFM.

The paper is organized as follows. In Section~\ref{sec2}, we state the
single-index model, discuss estimation of $g$ using
kernel estimating functions and of $\bolds{\beta}$ using profile
estimating functions, and investigate the problem of testing the
index using quasi-likelihood ratio. In Section~\ref{sec3} we provide a
computation algorithm for solving the estimating functions and
illustrate the method with simulation and practical studies. The
proofs are deferred to the \hyperref[app]{Appendix}.


\section{Estimating function method (EFM) and its large sample
properties}\label{sec2}

In this section, which is concerned with inference based on the
estimating function method, the model of interest is determined
through specification of mean and variance functions, up to an
unknown vector $\bolds{\beta}$ and an unknown function $g$. Except for
Gaussian data, model (\ref{model}) need not be a full
semiparametric likelihood specification. Note that the parameter
space $\Theta=\{{\bolds{\beta}}=(\beta_1,\ldots,\beta_d)^{\top} \dvtx
\|{\bolds{\beta}}\|=1, \beta_1>0, {\bolds{\beta}}\in{\R}^d\}$ means that
$\bolds{\beta}$ is on the boundary of a unit ball and it represents
therefore a nonregular problem. So we first choose an
identifiable parameterization which transforms the boundary of a
unit ball in ${\R}^d$ to the interior of a unit ball in
${\R}^{d-1}$. By eliminating $\beta_1$, the parameter space
$\Theta$ can be rearranged to a form $ \{((1-\sum_{r=2}^d
\beta_r^2)^{1/2},\beta_2,\ldots,\beta_d)^{\top}\dvtx
\sum_{r=2}^d \beta_r^2<1\}$. Then the derivatives of a
function with respect to
${\bolds{\beta}}^{(1)}=(\beta_2,\ldots,\beta_d)^{\top}$ are readily
obtained by chain rule and the classic assumptions on the
asymptotic normality hold after transformation. This
reparameterization is the key to analyzing the asymptotic
properties of the estimates for $\bolds{\beta}$ and to facilitating an
efficient computation algorithm. We will investigate the
estimation for $g$ and $\bolds{\beta}$ and propose a quasi-likelihood
method to test the statistical significance of certain variables
in the parametric component.

\subsection{The kernel estimating functions for the nonparametric part $g$}\label{sec2.1}

If $\bolds{\beta}$ is known, then we estimate $g(\cdot)$ and $g'(\cdot)$
using the local linear estimating functions. Let $h$ denote the
bandwidth parameter, and let $K(\cdot)$ denote the symmetric kernel
density function satisfying $K_h(\cdot)=h^{-1}K(\cdot/h)$. The
estimation method involves local linear approximation. Denote by
$\alpha_0$ and $\alpha_1$ the values of $g$ and $g'$ evaluating at
$t$, respectively. The local linear approximation for
$g({\bolds{\beta}}^{\top}{\mathbf{x}})$ in a neighborhood of $t$ is
$\tilde{g}({\bolds{\beta}}^{\top}{\mathbf{x}})=\alpha_0+\alpha_1({\bolds{\beta}}^{\top}{\mathbf{x}}-t)$. The estimators
$\hat{g}(t)$ and $\hat{g}'(t)$ are obtained by solving the kernel
estimating functions with respect to $\alpha_0, \alpha_1$:
%
\begin{equation}\label{eq-kernel}
\qquad \cases{
\displaystyle \sum_{j=1}^n K_h({\bolds{\beta}}^{\top}
{\mathbf{X}}_j-t) \mu'\{\tilde{g}({\bolds{\beta}}^{\top}{\mathbf{X}}_j)\}
V^{-1}\{\tilde{g}({\bolds{\beta}}^{\top}{\mathbf{X}}_j)\}\cr
\qquad {}\times[Y_j-\mu\{\tilde{g}({\bolds{\beta}}^{\top}{\mathbf{X}}_j)\}
]=0,\cr
\displaystyle \sum_{j=1}^n({\bolds{\beta}}^{\top}{\mathbf{X}}_j-t)K_h({\bolds{\beta}}^{\top}{\mathbf{X}}_j-t) \mu'\{\tilde{g}({\bolds{\beta}}^{\top}{\mathbf{X}}_j)\} V^{-1}\{\tilde{g}({\bolds{\beta}}^{\top}
{\mathbf{X}}_j)\}\cr
\qquad {}\times[Y_j-\mu\{\tilde{g}({\bolds{\beta}}^{\top}{\mathbf{X}}_j)\}]=0.
}
\end{equation}
Having estimated $\alpha_0, \alpha_1$ at $t$ as $\hat{\alpha}_0,
\hat{\alpha}_1$, the local linear estimators of $g(t)$ and $g'(t)$ are
$\hat{g}(t)=\hat{\alpha}_0$ and $\hat{g}'(t)=\hat{\alpha}_1$,
respectively.

The key to obtain the asymptotic normality of the estimates for
$\bolds{\beta}$ lies in the asymptotic properties of the estimated
nonparametric part. The following theorem will provide some useful
results. The following notation will be used. Let
$\mathcal{X}=\{{\mathbf{X}}_1,\ldots,{\mathbf{X}}_n\}$,
$\rho_l(z)=\{\mu'(z)\}^l V^{-1}(z)$ and ${\mathbf{J}}=\frac{\partial
{\bolds{\beta}}}{\partial{\bolds{\beta}}^{(1)}}$ the Jacobian matrix of size
$d\times(d-1)$ with
\[
{\mathbf{J}} = \pmatrix{
-{\bolds{\beta}}^{(1)\top}/\sqrt{1-\bigl\|{\bolds{\beta}}^{(1)}\bigr\|^2} \vspace*{2pt}\cr
{\mathbf{I}}_{d-1}
},\qquad
{\bolds{\beta}}^{(1)}=(\beta_2,\ldots,\beta_d)^{\top}.
\]
The moments of $K$ and $K^2$ are denoted, respectively, by,
$j=0,1,\ldots,$
\[
\gamma_j=\int t^j K(t)\,dt\quad  \mbox{and}\quad  \nu_j=\int t^j K^2(t)\,dt.
\]

\begin{prop}\label{prop-non}
Under regularity conditions \textup{(a), (b), (d)} and
\textup{(e)} given
in the \hyperref[app]{Appendix}, we have:
\begin{longlist}[(iii)]
\item[(i)] With $h\rightarrow0$,
$n\rightarrow\infty$ such that $h\rightarrow0$ and $nh\rightarrow
\infty$,
$\forall{\bolds{\beta}}\in\Theta$, the asymptotic conditional bias and
variance of
$\hat{g}$ are given by
%
\begin{eqnarray}
&&E \bigl\{\{\hat{g}({\bolds{\beta}}^{\top}{\mathbf{x}})-g({\bolds{\beta}}^{\top
}{\mathbf{x}})\}^2 |\mathcal{X} \bigr\}\nonumber \\
&&\qquad = \bigl\{\tfrac{1}{2}\gamma_2h^2g''({\bolds{\beta}}^{\top}{\mathbf{x}})
\bigr\}^2\nonumber\\[-8pt]\\[-8pt]
&&\qquad \quad {} + \nu_0\sigma^2/ 
[nhf_{{\bolds{\beta}}^{\top}{\mathbf{x}}}({\bolds{\beta}}^{\top} {\mathbf{x}})\rho
_2\{g({\bolds{\beta}}^{\top}{\mathbf{x}})\}] \nonumber \\ 
&&\qquad \quad {}+ {\mo}_P(h^4+n^{-1}h^{-1}).\nonumber
\end{eqnarray}
\item[(ii)] With $h\rightarrow0$, $n\rightarrow\infty$ such that
$h\rightarrow0$ and $nh^3\rightarrow\infty$, for the estimates of the
derivative $g'$, it holds that
%
\begin{eqnarray}
&&E \bigl\{\{\hat{g}'({\bolds{\beta}}^{\top}{\mathbf{x}})-g'({\bolds{\beta}}^{\top}
{\mathbf{x}})\}^2 |\mathcal{X} \bigr\}\nonumber\\
& &\qquad = \bigl\{\tfrac{1}{6}\gamma_4\gamma_2^{-1}h^2g'''({\bolds{\beta}}^{\top
}{\mathbf{x}})
\nonumber\\
&&\qquad \quad\hphantom{\bigl\{} {} +\tfrac{1}{2}(\gamma_4\gamma_2^{-1}-\gamma_2)h^2g''({\bolds{\beta}
}^{\top}
{\mathbf{x}})\nonumber\\[-8pt]\\[-8pt]\nonumber
&&\qquad \quad\hphantom{\bigl\{{}+{}} {} \times[ \rho_2'\{g({\bolds{\beta}}^{\top} {\mathbf{x}})\}/\rho_2\{g({\bolds{\beta}}^{\top} {\mathbf{x}})\}
 +f'_{{\bolds{\beta}}^{\top}{\mathbf{x}}}({\bolds{\beta}}^{\top}{\mathbf{x}})/f_{{\bolds{\beta}}^{\top}{\mathbf{x}}}({\bolds{\beta}}^{\top}{\mathbf{x}})] \bigr\}^2 \\ \nonumber
&& \qquad \quad {}+ \hfill\nu_2\gamma_2^{-2}\sigma^2/[nh^3f_{{\bolds{\beta}}^{\top
}{\mathbf{x}}}({\bolds{\beta}}^{\top}{\mathbf{x}})\rho_2\{g({\bolds{\beta}}^{\top}{\mathbf{x}})\}]\\ \nonumber
&&\qquad \quad {}+ {\mo}_P(h^4+n^{-1}h^{-3}).
\end{eqnarray}

%
\item[(iii)] With $h\rightarrow0$,
$n\rightarrow\infty$ such that $h\rightarrow0$ and $nh^3\rightarrow
\infty$,
we have that
%
\begin{equation}
 E \biggl\{\biggl\|\frac{\partial\hat{g}({\bolds{\beta}}^{\top}{\mathbf{x}})}{\partial
{\bolds{\beta}}^{(1)}}-g'({\bolds{\beta}}^{\top}{\mathbf{x}}){\mathbf{J}}^{\top} \{
{\mathbf{x}}-E({\mathbf{x}}|{\bolds{\beta}}^{\top}{\mathbf{x}})\}\biggr\|^2 \Big|\mathcal{X}\! \biggr\} =
\CO_P(h^4+n^{-1}h^{-3}).\hspace*{-35pt}
\end{equation}
\end{longlist}
\end{prop}

The proof of this proposition appears in the \hyperref[app]{Appendix}. Results (i) and
(ii) in
Proposition~\ref{prop-non} are routine and similar to
\citet{Carroll1998}. In the situation where $\sigma^2V=\sigma^2$
and the
function $\mu$ is identity, results (i) and (ii) coincide with those
given by
\citet{Fan1996}. From result (iii), it is seen that $\partial
\hat{g}({\bolds{\beta}}^{\top} {\mathbf{x}})/\partial{\bolds{\beta}}^{(1)}$
converges in
probability to $g'({\bolds{\beta}}^{\top}{\mathbf{x}}){\mathbf{J}}^{\top}\{{\mathbf{x}}-E({\mathbf{x}}|{\bolds{\beta}}^{\top} {\mathbf{x}})\}$, rather than $g'({\bolds{\beta}}^{\top
}{\mathbf{x}}){\mathbf{J}}^{\top}{\mathbf{x}}$ as if $g$ were known. That is,
$\lim_{n\rightarrow\infty}\{\partial\hat{g}({\bolds{\beta}}^{\top}
{\mathbf{x}})/\partial{\bolds{\beta}}^{(1)}\}\neq
\partial\{\lim_{n\rightarrow\infty}
\hat{g}({\bolds{\beta}}^{\top}{\mathbf{x}})\}/\partial{\bolds{\beta}}^{(1)}$,
which means
that the convergence in probability and the derivation of the sequence
$\hat{g}_n({\bolds{\beta}}^{\top} {\mathbf{x}})$ (as a function of $n$) cannot commute. This is primarily
caused by the fact that the support of the infinite-dimensional nuisance
parameter $g(\cdot)$ depends on the finite-dimensional projection parameter
$\bolds{\beta}$. In contrast, a semiparametric model where the support of the
nuisance parameter is independent of the finite-dimensional parameter
is a
partially linear regression model having form $Y={\mathbf{X}}^{\top}{\bolds{\theta}}+\eta(T)+\varepsilon$. It is easy to check
that the limit
of $\partial\hat{\eta}(T)/\partial{\bolds{\theta}}$ is equal to
$E({\mathbf{X}}|T)$,
which is the derivative of $ \lim_{n\rightarrow
\infty}\hat{\eta}(T)=E(Y|T)-E({\mathbf{X}}^{\top} |T){\bolds{\theta}}$ with
respect to
${\bolds{\theta}}$. Result (iii) ensures that the proposed estimator does not
require undersmoothing of $g(\cdot)$ to obtain a root-$n$ consistent estimator
for $\bolds{\beta}$ and it is also of its own interest in inference theory for
semiparametric models.

\subsection{\texorpdfstring{The asymptotic distribution for the estimates of the parametric part~$\beta$}
{The asymptotic distribution for the estimates of the parametric part beta}}\label{sec2.2}

We will now proceed to the estimation of $\bolds{\beta}\in\Theta$. We
need to
estimate the $(d-1)$-dimensional vector ${\bolds{\beta}}^{(1)}$, the
estimator of
which will be defined via
%
\begin{equation}\label{raw-eq}
\sum_{i=1}^n\bigl[\partial\mu\{\hat{g}({\bolds{\beta}}^{\top}{\mathbf{X}}_i)\}/\partial
{\bolds{\beta}}^{(1)}\bigr]{V^{-1}\{\hat{g}({\bolds{\beta}}^{\top}{\mathbf{X}}_i)\}}[Y_i-\mu\{\hat{g}({\bolds{\beta}}^{\top}{\mathbf{X}}_i)\}]=0.
\end{equation}
This is the direct analogue of the ``ideal'' estimating equation for
known~$g$, in that it is calculated by replacing $g(t)$ with
$\hat{g}(t)$. An asymptotically equivalent and easily computed
version of this equation is
%
\begin{eqnarray}\label{eq-profile}
 \hat{\mathbf{G}}({\bolds{\beta}})
&\stackrel{\mathrm{def}}=&\sum_{i=1}^n{\mathbf{J}}^{\top}
\hat{g}'({\bolds{\beta}}^{\top}{\mathbf{X}}_i)\{{\mathbf{X}}_i-\hat{\mathbf{h}}({\bolds{\beta}}^{\top}{\mathbf{X}}_i)\}\rho_1\{\hat{g}({\bolds{\beta}}^{\top
}{\mathbf{X}}_i)\}[Y_i-\mu\{\hat{g}({\bolds{\beta}}^{\top}{\mathbf{X}}_i)\}]\nonumber\hspace*{-35pt}\\[-8pt]\\[-8pt]
&=&0\nonumber\hspace*{-35pt}
\end{eqnarray}
with ${\mathbf{J}}=\frac{\partial{\bolds{\beta}}}{\partial{\bolds{\beta}
}^{(1)}}$ the
Jacobian mentioned above, $\hat{g}$ and $\hat{g}'$ are defined by\vspace*{-3pt}
(\ref{eq-kernel}), and $\hat{\mathbf{h}}(t)$ the local linear estimate for
${\mathbf{h}}(t)=E({\mathbf{X}}|{\bolds{\beta}}^{\top}{\mathbf{X}}=t)=(h_1(t),\ldots,\break
h_d(t))^{\top}$,
\[
\hat{\mathbf{h}}(t)=\sum_{i=1}^n b_{i}(t){\mathbf{X}}_{i} \bigg/\sum_{i=1}^n
b_i(t),
\]
where $b_i(t)= K_h({\bolds{\beta}}^{\top}{\mathbf{X}}_i-t)\{S_{n,2}(t)-({{\bolds{\beta}}^{\top}\mathbf{X}}_i-t)S_{n,1}(t)\}$,
$S_{n,k}=\break\sum_{i=1}^n K_h({\bolds{\beta}}^{\top}\times{\mathbf{X}}_i-t)({\bolds{\beta}}^{\top}{\mathbf{X}}_i-t)^k, k=1,2.$
We use (\ref{eq-profile}) to estimate ${\bolds{\beta}}^{(1)}$ in the
single-index model, and then use the fact that
$\beta_1=\sqrt{1-\|{\bolds{\beta}}^{(1)}\|^2}$ to obtain $\hat{\beta}_1$.
The use of (\ref{eq-profile}) constitutes in our view a new approach
to estimating single-index models; since (\ref{eq-profile})
involves smooth pilot estimation of $g$, $g'$ and ${\mathbf{h}}$ we call it the
Estimation Function Method (EFM) for ${\bolds{\beta}}$.

\begin{rem}
The estimating equations $\hat{\mathbf{G}}({\bolds{\beta}})$ can be represented as the gradient vector of the
following objective function:
\[
\hat{Q}({\bolds{\beta}})=\sum_{i=1}^n
Q[\mu\{\hat{g}({\bolds{\beta}}^{\top} {\mathbf{X}}_i)\},Y_i]
\]
with $Q[\mu,y]=\int_{\mu}^y \frac{s-y}{V\{\mu^{-1}(s)\}}\,ds$ and
$\mu^{-1}(\cdot)$ the inverse function of $\mu(\cdot)$. The\vspace*{-2pt}
existence of such a potential function makes $\hat{\mathbf{G}}({\bolds{\beta}})$ to inherit properties of the ideal likelihood score
function. Note that $\{\|{\bolds{\beta}}^{(1)}\|<1\}$ is an open,
connected subset of ${\R}^{d-1}$. By the regularity conditions
assumed on $\mu(\cdot), g(\cdot), V(\cdot)$ (for details see the
\hyperref[app]{Appendix}), we know that the quasi-likelihood function
$\hat{Q}({\bolds{\beta}})$ is twice continuously differentiable on
$\{\|{\bolds{\beta}}^{(1)}\|<1\}$ such that the global maximum of
$\hat{Q}({\bolds{\beta}})$ can be achieved at some point. One may
ask
whether the solution is unique and also consistent. Some elementary
calculations lead to the Hessian matrix
$\partial^2\hat{Q}({\bolds{\beta}}) /\partial{\bolds{\beta}}^{(1)}\partial
{\bolds{\beta}}^{(1)\top}$, because the partial derivative\vspace*{-3pt}
$\frac{\partial\mu\{\hat{g}({\bolds{\beta}}^{\top}{\mathbf{X}}_i)\}
}{\partial
{\bolds{\beta}}^{(1)}}=\mu'\{\hat{g}({\bolds{\beta}}^{\top}{\mathbf{X}}_i)\}\hat{g}'({\bolds{\beta}}^{\top}{\mathbf{X}}_i)\{{\mathbf{X}}_{i}
-\hat{\mathbf{h}}({\bolds{\beta}}^{\top}{\mathbf{X}}_i)\},$ then
\begin{eqnarray*}\label{Hmatrix}
&&
\frac{1}{n}\frac{\partial^2\hat{Q}({\bolds{\beta}
})}{\partial
{\bolds{\beta}}^{(1)}\partial{\bolds{\beta}}^{(1)\top}} \\
&&\quad =\frac{1}{n}\frac
{\partial\hat{\mathbf{G}}({\bolds{\beta}})}{\partial{\bolds{\beta}}^{(1)}}\\
&&\quad = \frac{1}{n}\sum_{i=1}^n \frac{\partial[{\mathbf{J}}^{\top}
\hat{g}'({\bolds{\beta}}^{\top}{\mathbf{X}}_i)\{{\mathbf{X}}_i-\hat{\mathbf{h}}({\bolds{\beta}}^{\top}{\mathbf{X}}_i)\}\rho_1\{\hat{g}({\bolds{\beta}}^{\top}{\mathbf{X}}_i)\}]}{\partial
{\bolds{\beta}}^{(1)}}[Y_i-\mu\{\hat{g}({\bolds{\beta}}^{\top}{\mathbf{X}}_i)\}]\\
&&\quad \quad {}- \frac{1}{n}\sum_{i=1}^n {\mathbf{J}}^{\top}
\hat{g}'({\bolds{\beta}}^{\top}{\mathbf{X}}_i)\{{\mathbf{X}}_i-\hat{\mathbf{h}}({\bolds{\beta}}^{\top}{\mathbf{X}}_i)\}\rho_1\{\hat{g}({\bolds{\beta}}^{\top}{\mathbf{X}}_i)\}\frac{\partial\mu\{\hat{g}({\bolds{\beta}}^{\top}{\mathbf{X}}_i)\}}{\partial{\bolds{\beta}}^{(1)}}\\
&&\quad =\frac{1}{n}\sum_{i=1}^n \biggl[-\frac{\partial\{{\bolds{\beta}
}^{(1)}/\sqrt{1-\|{\bolds{\beta}}^{(1)}\|^2}\}}
{\partial{\bolds{\beta}}^{(1)}} \hat{g}'({\bolds{\beta}}^{\top}{\mathbf{X}}_i)\{{\mathbf{X}}_{1i}-\hat{ h}_1({\bolds{\beta}}^{\top}{\mathbf{X}}_i)\}\rho_1\{\hat{g}({\bolds{\beta}}^{\top}{\mathbf{X}}_i)\}\\
&&\quad \quad\hphantom{\frac{1}{n}\sum_{i=1}^n \biggl[} {}+{\mathbf{J}}^{\top}\{{\mathbf{X}}_{i}-\hat{\mathbf{h}}({\bolds{\beta}}^{\top}{\mathbf{X}}_i)\}\frac{\partial\hat{g}'({\bolds{\beta}}^{\top}{\mathbf{X}}_i)}
{\partial
{\bolds{\beta}}^{(1)\top}}\rho_1\{\hat{g}({\bolds{\beta}}^{\top}{\mathbf{X}}_i)\}\\
&&\quad \quad\hphantom{\frac{1}{n}\sum_{i=1}^n \biggl[} {}+{\mathbf{J}}^{\top}\hat{g}'({\bolds{\beta}}^{\top}{\mathbf{X}}_i)\{{\mathbf{X}}_{i}-\hat{\mathbf{h}}({\bolds{\beta}}^{\top}{\mathbf{X}}_i)\}\frac{\partial\rho_1\{\hat{g}({\bolds{\beta}}^{\top}{\mathbf{X}}_i)\}}
{\partial{\bolds{\beta}}^{(1)\top}}\\
&&\quad \quad\hspace*{143pt} {} -{\mathbf{J}}^{\top}\hat{g}'({\bolds{\beta}}^{\top}{\mathbf{X}}_i)
\frac{\partial\hat{\mathbf{h}}({\bolds{\beta}}^{\top}{\mathbf{X}}_i)}{\partial
{\bolds{\beta}}^{(1)}}\rho_1\{\hat{g}({\bolds{\beta}}^{\top}{\mathbf{X}}_i)\} \biggr] \\
&&\quad \quad\hphantom{\frac{1}{n}\sum_{i=1}^n}
{}\times[Y_i-\mu\{\hat{g}({\bolds{\beta}}^{\top}{\mathbf{X}}_i)\}]\\
&&\quad \quad {} - \frac{1}{n}\sum_{i=1}^n {\mathbf{J}}^{\top}
\hat{g}'^2({\bolds{\beta}}^{\top}{\mathbf{X}}_i)\{{\mathbf{X}}_i-\hat{\mathbf{h}}({\bolds{\beta}}^{\top}{\mathbf{X}}_i)\}\{{\mathbf{X}}_i-\hat{\mathbf{h}}({\bolds{\beta}}^{\top}{\mathbf{X}}_i)\}^{\top}\rho_2\{\hat{g}({\bolds{\beta}}^{\top}{\mathbf{X}}_i)\}{\mathbf{J}}.
\end{eqnarray*}
By the regularity conditions in the \hyperref[app]{Appendix}, the multipliers of
the residuals $[Y_i-\mu\{\hat{g}({\bolds{\beta}}^{\top}{\mathbf{X}}_i)\}]$ in
the first sum of (\ref{Hmatrix}) are bounded. Mimicking the proof of
Proposition~\ref{prop-non}, the first sum can be shown to converge to 0 in
probability as $n$ goes to infinity. The second sum converges to a
negative semidefinite matrix. If the Hessian matrix
$\frac{1}{n}\frac{\partial^2\hat{Q}({\bolds{\beta}})}{\partial
{\bolds{\beta}}^{(1)}\partial{\bolds{\beta}}^{(1)\top}}$ is negative definite
for all values of ${\bolds{\beta}}^{(1)}$, $\hat{\mathbf{G}}({\bolds{\beta}})$ has
a unique root.
At sample level, however, estimating functions may have more than
one root. For the EFM method, the quasi-likelihood
$\hat{Q}({\bolds{\beta}})$ exists, which can be used to distinguish local
maxima from minima. Thus, we suppose (\ref{eq-profile}) has a unique
solution in the following context.
\end{rem}

\begin{rem}
It can be seen from the proof in the \hyperref[app]{Appendix} that the population
version of $\hat{\mathbf{G}}({\bolds{\beta}})$ is
%
\begin{equation}\label{peqp}
 {\mathbf{G}}({\bolds{\beta}})\!=\! \sum_{i=1}^n\!{\mathbf{J}}^{\top}
g'({\bolds{\beta}}^{\top}{\mathbf{X}}_i)\{{\mathbf{X}}_i\!-\!{\mathbf{h}}
({\bolds{\beta}}^{\top}{\mathbf{X}}_i)\}\rho_1\{g({\bolds{\beta}}^{\top}{\mathbf{X}}_i)\}
[Y_i\!-\!\mu\{g({\bolds{\beta}}^{\top}{\mathbf{X}}_i)\}],\hspace*{-35pt}
\end{equation}
which is obtained by replacing $\hat{g}, \hat{g}', \hat{\mathbf{h}}$ with
$g, g', {\mathbf{h}}$ in (\ref{eq-profile}). One important property of
(\ref{peqp}) is that the second Bartlett identity holds, for any
$\bolds{\beta}$:
\[
E\{{\mathbf{G}}({\bolds{\beta}}){\mathbf{G}}^{\top}
({\bolds{\beta}})\}=-E\biggl\{\frac{\partial{{\mathbf{G}}({\bolds{\beta}})}}{\partial{{\bolds{\beta}}^{(1)}}}\biggr\}.
\]
This property makes the semiparametric efficiency of
the EFM (\ref{eq-profile}) possible.
\end{rem}


Let ${\bolds{\beta}}^{0}=(\beta_1^0,
{{\bolds{\beta}}^{(1)0}}^{\top})^{\top}$ denote the true parameter and
${\mathbf{B}}^+$ denote the Moore--Penrose inverse of any given matrix
${\mathbf{B}}$. We have the following asymptotic result for the estimator
$\hat{\bolds{\beta}}{}^{(1)}$.

\begin{thm}\label{asymptotic-normality} Assume the
estimating function (\ref{eq-profile}) has\vspace*{-1pt} a unique solution and
denote it by $\hat{\bolds{\beta}}{}^{(1)}$. If the regularity conditions
\textup{(a)--(e)} in the \hyperref[app]{Appendix} are satisfied, the
following results hold:
\begin{longlist}[(ii)]
\item[(i)]
With $h\rightarrow0$, $n\rightarrow\infty$ such that
$(nh)^{-1}\log(1/h)\rightarrow0$, $\hat{\bolds{\beta}}{}^{(1)}$ converges
in probability to the true
parameter ${\bolds{\beta}}^{(1)0}$.
\item[(ii)] If $nh^6\rightarrow0$ and $nh^4\rightarrow\infty$,
%
\begin{equation}
\sqrt{n}\bigl(\hat{\bolds{\beta}}{}^{(1)}-{\bolds{\beta}}^{(1)0}\bigr)\stackrel
{\mathcal{L}}
\longrightarrow N_{d-1}({\bolds{0}},{\bolds{\Sigma}}_{{\bolds{\beta}}^{(1)0}}),\vadjust{\goodbreak}
\end{equation}
where ${\bolds{\Sigma}}_{{\bolds{\beta}}^{(1)0}}=\{{\mathbf{J}}^{\top} {\bolds{\Omega}}{\mathbf{J}}\}^{+}|_{{\bolds{\beta}}^{(1)}={\bolds{\beta}}^{(1)0}}$, ${\mathbf{J}}=\frac{\partial{\bolds{\beta}}}{\partial{\bolds{\beta}}^{(1)}}$ and
\[
{\bolds{\Omega}} =
E [\{{\mathbf{XX}}^{\top}-E({\mathbf{X}}|{\bolds{\beta}}^{\top}{\mathbf{X}})E({\mathbf{X}}^{\top} |{\bolds{\beta}}^{\top}{\mathbf{X}})\}\rho_2\{g({\bolds{\beta}}^{\top}{\mathbf{X}})\}\{g'({\bolds{\beta}}^{\top
}{\mathbf{X}})\}^2/\sigma^2 ].
\]
\end{longlist}
\end{thm}

\begin{rem}
Note that ${\bolds{\beta}}^{\top}{\bolds{\Omega}}{\bolds{\beta}}=0$, so the
nonnegative matrix
${\bolds{\Omega}}$ degenerates in the direction of $\bolds{\beta}$. If the
mean function
$\mu$ is the identity function and the variance function is equal to a scale
constant, that is, $\mu\{g({\bolds{\beta}}^{\top}{\mathbf{X}})\}=g({\bolds{\beta}
}^{\top} {\mathbf{X}})$, $\sigma^2V\{g({\bolds{\beta}}^{\top}{\mathbf{X}})\}=\sigma^2$, the
matrix ${\bolds{\Omega}}$ in Theorem~\ref{asymptotic-normality} reduces to be
\[
{\bolds{\Omega}} = E [\{{\mathbf{X}}{\mathbf{X}}^{\top} -{E}({\mathbf{X}}|{\bolds{\beta}
}^{\top}{\mathbf{X}}){E}({\mathbf{X}}^{\top}|{\bolds{\beta}}^{\top} {\mathbf{X}})\} \{g'({\bolds{\beta}
}^{\top} {\mathbf{X}})\}^2/\sigma^2 ].
\]
\end{rem}

Technically speaking, Theorem~\ref{asymptotic-normality} shows that an undersmoothing approach is
unnecessary and that root-$n$ consistency can be achieved. The asymptotic
covariance ${\bolds{\Sigma}}_{{\bolds{\beta}}^{(1)0}}$ in general can be
estimated by\vspace*{-1.5pt}
replacing terms in its expression by estimates of those terms. The asymptotic
normality of $\hat{\bolds{\beta}}=(\hat{\beta}_1, {\hat{\bolds{\beta}
}{}^{(1)}}^{\top})^{\top} $
will follow from Theorem~\ref{asymptotic-normality} with a simple application
of the multivariate delta-method, since
$\hat{\beta}_1=\sqrt{1-\|\hat{\bolds{\beta}}{}^{(1)}\|^2}$. According to
the results
of \citet{Carroll1997}, the asymptotic variance of their estimator is
${\bolds{\Omega}}^+$. Define the block partition of matrix ${\bolds{\Omega}}$ as follows:
%
\begin{equation}
{\bolds{\Omega}}= \pmatrix{
{\bolds{\Omega}}_{11} & {\bolds{\Omega}}_{12} \cr
{\bolds{\Omega}}_{21} & {\bolds{\Omega}}_{22}
},
\end{equation}
where ${\bolds{\Omega}}_{11}$ is a positive constant,
${\bolds{\Omega}}_{12}$ is a $(d-1)$-dimensional row vector,
${\bolds{\Omega}}_{21}$ is a $(d-1)$-dimensional column vector and
${\bolds{\Omega}}_{22}$ is a $(d-1)\times(d-1)$ nonnegative definite
matrix.

\begin{cor}\label{variance-comparision}
Under the conditions of Theorem
\ref{asymptotic-normality}, we have
%
\begin{equation}\label{asym-var}
\sqrt{n}(\hat{\bolds{\beta}}-{\bolds{\beta}}^{0})\stackrel{\mathcal{L}}
\longrightarrow
N_p({\bolds{0}},{\bolds{\Sigma}}_{{\bolds{\beta}}^{0}})
\end{equation}
with ${\bolds{\Sigma}}_{{\bolds{\beta}}^{0}}={\mathbf{J}}\{{\mathbf{J}}^{\top}
{\bolds{\Omega}}{\mathbf{J}}\}^{+}{\mathbf{J}}^{\top} |_{{\bolds{\beta}}={\bolds{\beta}}^{0}}$.
Further,
\[
{\bolds{\Sigma}}_{{\bolds{\beta}}^{0}}\leq\bolds{\Omega}^+ |_{{\bolds{\beta}}={\bolds{\beta}}^{0}}
\]
and a strict less-than sign holds when $\det({\bolds{\Omega}}_{22})=0$.
That is, in this case EFM is more efficient than that of Carroll et
al. (\citeyear{Carroll1997}).
\end{cor}

The possible smaller limiting variance derived from the
EFM approach
partly benefits from the reparameterization so that the
quasi-likelihood can be adopted. As we know, the
quasi-likelihood is often of optimal property. In
contrast, most existing methods treat the estimation of
${\bolds{\beta}}$ as if it were
done in the framework of linear dimension reduction. The
target of
linear dimension reduction is to find the directions
that can
linearly transform the original variables vector into a
vector of
one less dimension. For example, ADE and SIR are two
relevant
methods. However, when the link function $\mu(\cdot)$
is identity,
the limiting variance derived here\vadjust{\goodbreak} may not be smaller or
equal to
the ones of \citet{WangXueZhuChong2010} and \citet{Chang2010}
when the quasi-likelihood of (2.5) is applied.

\subsection{Profile quasi-likelihood ratio test}\label{sec2.3}

In applications, it is important to test the statistical significance
of added
predictors in a regression model. Here we establish a quasi-likelihood ratio
statistic to test the significance of certain variables in the linear index.
The null hypothesis that the model is correct is tested against a full model
alternative. \citet{FanJiang2007} gave a recent review about generalized
likelihood ratio tests. Bootstrap tests for nonparametric regression,
generalized partially linear models and single-index models have been
systematically investigated [see \citet{HardleMammen1993}, \citet{HardleMammenMuller1998}, \citet{HardleMammenProenca2001}].
Consider the testing problem:
%
\begin{eqnarray}
&&H_0\dvtx g(\cdot)=g\Biggl(\sum_{k=1}^r\beta_kX_k\Biggr)\nonumber\\[-8pt]\\[-8pt]
&&\qquad  \longleftrightarrow\quad
H_1\dvtx g(\cdot)=g\Biggl(\sum_{k=1}^r\beta_kX_k+\sum_{k=r+1}^d\beta_kX_k\Biggr).\nonumber
\end{eqnarray}
We mainly focus on testing $\beta_k=0, k=r+1,\ldots,d$, though the
following test procedure can be easily extended to a general linear testing
${\mathbf{B}}\tilde{\bolds{\beta}}=0$ where ${\mathbf{B}}$ is a known matrix with
full row
rank and $\tilde{\bolds{\beta}}=(\beta_{r+1},\ldots,\beta_{d})^{\top}
$. The profile
quasi-likelihood ratio test is defined by
%
\begin{equation}
T_n=2\Bigl\{\sup_{{\bolds{\beta}}\in\Theta}\hat{Q}({\bolds{\beta}}) -
\sup_{{\bolds{\beta}}\in\Theta,\widetilde{\bolds{\beta}}=0}\hat{Q}(\bolds{\beta})\Bigr\},
\end{equation}
where $\hat{Q}({\bolds{\beta}})=\sum_{i=1}^n
Q[\mu\{\hat{g}({\bolds{\beta}}^{\top} {\mathbf{X}}_i)\},Y_i],
Q[\mu,y]=\int_{\mu}^y\frac{s-y}{V\{\mu^{-1}(s)\}}\,ds$ and
$\mu^{-1}(\cdot)$ is the inverse function of $\mu(\cdot)$. The
following Wilks type theorem shows that the distribution of $T_n$
is asymptotically chi-squared and independent of nuisance
parameters.

\begin{thm}\label{test}
Under the assumptions of Theorem~\ref{asymptotic-normality},
if $\beta_k=0, k=r+1,\ldots,d$, then
%
\begin{equation}
T_n\stackrel{\mathcal{L}} \longrightarrow\chi^2(d-r).
\end{equation}
\end{thm}


\section{Numerical studies}\label{sec3}

\subsection{Computation of the estimates}\label{sec3.1}

Solving the joint estimating\vspace*{1pt} equations (\ref{eq-kernel}) and (\ref
{eq-profile})
poses some interesting challenges, since the functions\vspace*{1pt}
$\hat{g}({\bolds{\beta}}^{\top} {\mathbf{X}})$ and $\hat{g}'({\bolds{\beta}
}^{\top}{\mathbf{X}})$
depend on $\bolds{\beta}$ implicitly. Treating ${\bolds{\beta}}^{\top}X$ as
a new
predictor (with given $\bolds{\beta}$), (\ref{eq-kernel}) gives us
$\hat{g}, \hat{g}'$ as in \citet{Fan1995}. We therefore focus on\vadjust{\goodbreak}
(\ref{eq-profile}), as estimating equations. It cannot be solved explicitly,
and hence one needs to find solutions using numerical methods. The
Newton--Raphson algorithm is one of the popular and successful methods for
finding roots. However, the computational speed of this algorithm crucially
depends on the initial value. We propose therefore a
fixed-point iterative
algorithm that is not very sensitive to starting values and is adaptive to
larger dimension.
It is worth noting that this algorithm can be implemented in the case
that $d$ is slightly larger than $n$, because the
resultant procedure only involves one-dimensional
nonparametric smoothers, thereby avoiding the data sparsity problem
caused by
high dimensionality.

Rewrite the estimating functions as $\hat{\mathbf{G}}({\bolds{\beta}})={\mathbf{J}}^{\top} \hat{\mathbf{F}}({\bolds{\beta}})$ with
\[
\hat{\mathbf{F}}({\bolds{\beta}})=(\hat{F}_1({\bolds{\beta}}),\ldots,\hat{F}_d({\bolds{\beta}}))^{\top}
\]
and
\begin{eqnarray*}
\hat{F}_s({\bolds{\beta}})&=&\sum_{i=1}^n
\{X_{si}-\hat{h}_s({\bolds{\beta}}^{\top} {\mathbf{X}}_{i})\}\mu'\{\hat{g}({\bolds{\beta}}^{\top} {\mathbf{X}}_i)\}
\hat{g}'({\bolds{\beta}}^{\top} {\mathbf{X}}_i)V^{-1}\{\hat{g}({\bolds{\beta}}^{\top} {\mathbf{X}}_i)\}\\
&&{}\times[Y_i-\mu\{\hat{g}({\bolds{\beta}}^{\top} {\mathbf{X}}_i)\}].
\end{eqnarray*}
Setting
$\hat{\mathbf{G}}({\bolds{\beta}})=0$, we have that
%
\begin{equation}
\cases{
-\beta_2\hat{F}_1({\bolds{\beta}})/\sqrt{1-\bigl\|{\bolds{\beta}}^{(1)}\bigr\|
^2}+\hat{F}_2({\bolds{\beta}})=0,\cr
-\beta_3\hat{F}_1({\bolds{\beta}})/\sqrt{1-\bigl\|{\bolds{\beta}}^{(1)}\bigr\|
^2}+\hat{F}_3({\bolds{\beta}})=0,\cr
\cdots\cr
-\beta_d\hat{F}_1({\bolds{\beta}})/\sqrt{1-\bigl\|{\bolds{\beta}}^{(1)}\bigr\|
^2}+\hat{F}_d({\bolds{\beta}})=0.
}
\end{equation}
Note that
$\|{\bolds{\beta}}^{(1)}\|^2=\sum_{r=2}^d\beta_r^2$,
$\beta_1=\sqrt{1-\|{\bolds{\beta}}^{(1)}\|^2}$ and after some simple
calculations, we can get that
\[
\cases{
\beta_1=|\hat{F}_1({\bolds{\beta}})|
/\|\hat{\mathbf{F}}({\bolds{\beta}})\|,&\quad $ s=1$,\vspace*{2pt}\cr
\beta_s^2=\hat{F}_s^2({\bolds{\beta}})/\|\hat{\mathbf{F}}({\bolds{\beta}})\|^2,&\quad
$s\geq2$,
}
\]
and $\operatorname{sign}\{\beta_s\hat{F}_1({\bolds{\beta}})\}=\operatorname{sign}\{\hat{F}_s({\bolds{\beta}})\}, s\geq2$. The above equation can
also be rewritten as
%
\begin{equation}\label{iter}
{\bolds{\beta}} \frac{\hat{F}_1({\bolds{\beta}}) }{\|\hat{\mathbf{F}}({\bolds{\beta}})\|} = \frac{|\hat{F}_1({\bolds{\beta}})| }{\|\hat{\mathbf{F}}({\bolds{\beta}})\|}\times\frac{\hat{\mathbf{F}}({\bolds{\beta}})}{\|\hat
{\mathbf{F}}({\bolds{\beta}})\|}.
\end{equation}
Then solving the equation (\ref{eq-profile}) is equivalent to
finding a fixed point for (\ref{iter}). Though
$\|{\bolds{\beta}}^{(1)}\|<1$ holds almost surely in (\ref{iter}) and
always $\|{\bolds{\beta}}\|=1$, there will be some trouble if
(\ref{iter}) is directly used as iterative equations. Note that the
value of $\|\hat{\mathbf{F}}({\bolds{\beta}})\|$ is used as denominator that
may sometimes be small, which potentially makes the algorithm
unstable. On the other hand, the convergence rate of the fixed-point
iterative algorithm derived from (\ref{iter}) depends on $L$, where
$ \|\frac{\partial\{\hat{\mathbf{F}}({\bolds{\beta}})|/\|\hat{\mathbf{F}}({\bolds{\beta}})\|\}}{\partial{\bolds{\beta}}} \|\leq L$. For a fast
convergence rate, it technically needs a shrinkage value $L$. An ad
hoc fix introduces a constant $M$, adding $M {\bolds{\beta}}$ on both
sides of (\ref{iter}) and dividing by $\hat{F}_1({\bolds{\beta}})
/\|\hat{\mathbf{F}}({\bolds{\beta}})\|+M$:
\[
{\bolds{\beta}} = \frac{M}{\hat{F}_1({\bolds{\beta}})/\|\hat{\mathbf{F}}({\bolds{\beta}})\|+M}{\bolds{\beta}}+
\frac{|\hat{F}_1({\bolds{\beta}})|/\|\hat{\mathbf{F}}({\bolds{\beta}})\|^2}{\hat{F}_1({\bolds{\beta}}) /\|\hat{\mathbf{F}}({\bolds{\beta}})\|+M} \hat{\mathbf{F}}({\bolds{\beta}}),
\]
where $M$ is chosen such that $\hat{F}_1({\bolds{\beta}}) /\|\hat{\mathbf{F}}({\bolds{\beta}})\|+M \neq0$. In addition, to accelerate the rate of
convergence, we reduce the derivative of the term on the right-hand
side of the above equality, which can be achieved by choosing some
appropriate~$M$. This is the iteration formulation in Step 2. Here
the norm of ${\bolds{\beta}}_{\mathit{new}}$ is not equal to~1 and we have to
normalize it again. Since the iteration in Step 2 makes
${\bolds{\beta}}_{\mathit{new}}$ to violate the identifiability
constraint with norm 1, we design (\ref{iter}) to include the whole
${\bolds{\beta}}$
vector. The possibility of renormalization for ${\bolds{\beta}}_{\mathit{new}}$
avoids the difficulty of controlling $\|{\bolds{\beta}}_{\mathit{new}}^{(1)}\|<1$
in each iteration in Step 2.

Based on these observations, the fixed-point iterative algorithm is
summarized as:

\begin{longlist}[\textit{Step 0}.]
\item[\textit{Step 0}.] Choose initial values for $\bolds{\beta}$, denoted by ${\bolds{\beta}}_{\mathit{old}}$.

\item[\textit{Step 1}.] Solve the estimating equation (\ref{eq-kernel}) with
respect to
${\bolds{\alpha}}$, which yields $\hat{g}({\bolds{\beta}}_{\mathit{old}}^{\top}
{\mathbf{x}}_i)$ and $\hat{g}'({\bolds{\beta}}_{\mathit{old}}^{\top}{\mathbf{x}}_i)$,
$1\leq i\leq n$.

\item[\textit{Step 2}.] Update ${\bolds{\beta}}_{\mathit{old}}$ with ${\bolds{\beta}}_{\mathit{old}}={\bolds{\beta}}_{\mathit{new}}/\|{\bolds{\beta}}_{\mathit{new}}\|$ by solving
the equation (\ref{eq-profile}) in the fixed-point iteration
\[
{\bolds{\beta}}_{\mathit{new}} =
\frac{M}{\hat{F}_1({\bolds{\beta}}_{\mathit{old}})/\|\hat{F}({\bolds{\beta}
}_{\mathit{old}})\|+M}{\bolds{\beta}}_{\mathit{old}}+
\frac{|\hat{F}_1({\bolds{\beta}}_{\mathit{old}})|/\|\hat{F}({\bolds{\beta}}_{\mathit{old}})\|
^2}{\hat{F}_1({\bolds{\beta}}_{\mathit{old}})
/\|\hat{F}({\bolds{\beta}}_{\mathit{old}})\|+M} \hat{\mathbf{F}}({\bolds{\beta}}_{\mathit{old}}),
\]
where $M$ is a constant satisfying $\hat{F}_1({\bolds{\beta}})
/\|\hat{F}({\bolds{\beta}})\|+M \neq0$ for any $\bolds{\beta}$.

\item[\textit{Step 3}.] Repeat Steps 1 and 2 until $\max_{1\leq s \leq d}|\beta
_{\mathit{new},s}-\beta_{\mathit{old},s}|\leq{\mathit{tol}}$
is met with ${\mathit{tol}}$ being a prescribed tolerance.
\end{longlist}

The final vector ${\bolds{\beta}}_{\mathit{new}}/\|{\bolds{\beta}}_{\mathit{new}}\|$ is the
estimator of ${\bolds{\beta}}^0$. Similarly to other direct estimation
methods \citep{Horowitz1996}, the preceding calculation is easy to
implement. Empirically the initial value for ${\bolds{\beta}}$,
$(1,1,\ldots,1)^{\top} /\break\sqrt{d}$ can be used in the calculations.
The Epanechnikov kernel function $K(t)=3/4(1-t^2)I(|t|\leq1)$ is
used. The bandwidth involved in Step 1 can be chosen to be optimal
for estimation of $\hat{g}(t)$ and $\hat{g}'(t)$ based on the
observations $\{{\bolds{\beta}}_{\mathit{old}}^{\top} {\mathbf{X}}_i, Y_i\}$. So the
standard bandwidth selection methods, such as $K$-fold
cross-validation, generalized cross-validation (GCV) and the rule of
thumb, can be adopted. In this step, we recommend $K$-fold
cross-validation to determine the optimal bandwidth using the
quasi-likelihood as a criterion function. The $K$-fold cross-validation
is not too computationally intensive while making $K$ not
take too large values (e.g., $K=5$). Here we recommend trying
a number of smoothing parameters that smooth the data and picking the
one that seems most reasonable. As an adjustment
factor, $M$ will increase the stability of iteration. Ideally, in
each iteration an optimum value for $M$ should be chosen
guaranteeing that the derivative on the right-hand side of the
iteration formulation in Step 2 is close to zero. Following this
idea, $M$ will be depending the changes of ${\bolds{\beta}}$ and
$\hat{\mathbf{F}}({\bolds{\beta}}) /\|\hat{\mathbf{F}}({\bolds{\beta}})\|$. This will be
an expensive task due to the computation for the derivative on the
right-hand side of the iteration formulation in Step 2. We therefore
consider $M$ as constant nonvarying in each iteration, and select
$M$ by the $K$-fold cross-validation method, according to minimizing
the model prediction error. When the dimension $d$ gets larger, $M$
will get smaller. In our simulation runs, we empirically search $M$
in the interval $[2/\sqrt{d},d/2]$. This choice gives pretty good
practical performance.

\subsection{Simulation results}\label{sec3.2}

\begin{exm}[(Continuous response)]\label{ex1}
We report a
simulation study to investigate the finite-sample performance of the
proposed estimator and compare it with the rMAVE [refined MAVE; for
details see \citet{Xia2002}] estimator and the EDR estimator
[see \citet{Hristacheb2001}, \citet{Polzehl2009}]. We consider the
following model similar to that used in \citet{Xia2006}:
%
\begin{eqnarray}\label{simulation-model}
E(Y|{\bolds{\beta}}^{\top} {\mathbf{X}})&=&g({\bolds{\beta}}^{\top}{\mathbf{X}}),\qquad
g({\bolds{\beta}}^{\top} {\mathbf{X}}) = ({\bolds{\beta}}^{\top} {\mathbf{X}})^2
\exp({\bolds{\beta}}^{\top}{\mathbf{X}});\nonumber\\[-8pt]\\[-8pt]
\operatorname{Var}(Y|{\bolds{\beta}}^{\top} {\mathbf{X}})&=&\sigma^2, \qquad  \sigma=0.1.\nonumber
\end{eqnarray}
Let the true parameter ${\bolds{\beta}}=(2,1,0,\ldots,0)^{\top}
/\sqrt{5}$. Two sets of designs for $\mathbf{X}$ are considered: Design
(A) and Design (B). In Design (A), $(X_s+1)/2\sim \operatorname{Beta}(\tau,1)$,
$1\leq s \leq d$ and, in Design (B), $(X_1+1)/2\sim \operatorname{Beta}(\tau,1)$ and
$P(X_s=\pm0.5)=0.5$, $s=2,3,4,\ldots,d$. The data generated in Design
(A) are not elliptically symmetric. All the components of Design (B)
are discrete except for the first component $X_1$. $Y$ is generated
\begin{table}
\tablewidth=270pt
\caption{Average estimation errors $\sum_{s=1}^d
|\hat{\beta}_s-\beta_s|$ for model (\protect\ref{simulation-model})}\label{simulation1}
\begin{tabular*}{\tablewidth}{@{\extracolsep{\fill}}lccccccc@{}}
\hline
 &&\multicolumn{3}{c}{\textbf{Design (A)}} & \multicolumn{3}{c@{}}{\textbf{Design (B)}}\\[-5pt]
 &&\multicolumn{3}{c}{\hrulefill} &\multicolumn{3}{c@{}}{\hrulefill}\\
$\bolds{d}$ & $\bolds{\tau}$ & \textbf{rMAVE} & \textbf{EDR} & \textbf{EFM}& \textbf{rMAVE} & \textbf{EDR} & \textbf{EFM}
\\
\hline
10 & 0.75& 0.0559\tabnoteref{table1} & 0.0520& 0.0792 & 0.0522\tabnoteref{table1}&0.0662 & 0.0690 \\
10 & 1.5\hphantom{0} & 0.0323\tabnoteref{table1} & 0.0316& 0.0298 & 0.0417\tabnoteref{table1}& 0.0593 & 0.0457\\
50 & 0.75& 0.9900\hphantom{\tabnoteref{table1}} & 0.7271& 0.5425 & 0.9780\hphantom{\tabnoteref{table1}} & 0.7712 & 0.4515 \\
50 & 1.5\hphantom{0} & 0.3776\hphantom{\tabnoteref{table1}} & 0.3062& 0.1796 & 0.4693\hphantom{\tabnoteref{table1}} & 0.4103 & 0.2211 \\
\hline
\end{tabular*}
\tabnotetext{table1}{The values
are adopted from \citet{Xia2006}.}
\end{table}%
from a normal distribution. This simulation data set consists of
$400$ observations with $250$ replications. The results are shown in
Table~\ref{simulation1}. All rMAVE, EDR and EFM estimates are close
to the true parameter vector for $d=10$. However, the average
estimation errors from rMAVE and EDR estimates for $d=50$ are about
$2$ and $1.5$ times as large as those of the EFM estimates,
respectively. This indicates that the fixed-point algorithm is more
adaptive to high dimension.
\end{exm}

\begin{exm}[(Binary response)]\label{ex2}
This simulation design assumes an
underlying single-index model for binary responses with
%
\begin{eqnarray}\label{simulation-GSIM}
P(Y=1|{\mathbf{X}})&=&\mu\{g({\bolds{\beta}}^{\top} {\mathbf{X}})\}= \exp\{g({\bolds{\beta}}^{\top}
{\mathbf{X}})\}/[1+\exp\{g({\bolds{\beta}}^{\top}{\mathbf{X}})\}],\nonumber\\[-8pt]\\[-8pt]
g({\bolds{\beta}}^{\top} {\mathbf{X}})&=&\exp(5{\bolds{\beta}}^{\top}{\mathbf{X}}-2)/
\{1+\exp(5{\bolds{\beta}}^{\top} {\mathbf{X}}-3)\}-1.5.\nonumber
\end{eqnarray}
The underlying coefficients are assumed to be
${\bolds{\beta}}=(2,1,0,\ldots,0)^{\top}/\sqrt{5}$. We consider two sets
of designs: Design (C) and Design (D). In Design (C), $X_1$ and
$X_2$ follow the uniform distribution $U(-2,2)$. In Design (D),
$X_1$ is also assumed to be uniformly distributed in interval
$(-2,2)$ and $(X_2+1)/2\sim \operatorname{Beta}(1,1)$. Similar designs for
generalized partially linear single-index models are assumed in
\citet{KaneHoltAllen2004}. Here a sample size of $700$ is used
for the case $d=10$ and 3,000 is used for $d=50$. Different sample
sizes from Example~\ref{ex1} are used due to varying complexity of the two
examples. For this example, $250$ replications are simulated and the
results are displayed in Table~\ref{simulation2}. In this set of
simulations, the average estimation errors from rMAVE estimates and
EDR estimates are about $1.5$ and $1.2$ times as large as EFM
estimates, under both Design (C) and Design (D) for $d=10$ or
\begin{table}[b]
\caption{Average estimation errors $\sum_{s=1}^d
|\hat{\beta}_s-\beta_s|$ for model (\protect\ref{simulation-GSIM})}\label{simulation2}
\begin{tabular}{@{}lcccccc@{}}
\hline
& \multicolumn{3}{c}{\textbf{Design (C)}} & \multicolumn{3}{c@{}}{\textbf{Design (D)}}\\[-5pt]
& \multicolumn{3}{c}{\hrulefill}& \multicolumn{3}{c@{}}{\hrulefill}\\
$\bolds{d}$ & \textbf{rMAVE} & \textbf{EDR} & \textbf{EFM} &\textbf{rMAVE} & \textbf{EDR} &
\textbf{EFM}\\
\hline
10 & 0.5017 & 0.5281& 0.4564 & 0.9614 & 0.9574 & 0.7415 \\
50 & 2.0991 & 1.2695& 1.1744 & 2.5040 & 2.4846 & 1.9908 \\
\hline
\end{tabular}
\end{table}%
$d=50$. The values in the row marked by $d=50$ look a little bigger.
However, it is reasonable because the number of summands in the
average estimate error for $d=50$ is five times as large as that for
$d=10$. Again it appears that the EFM procedure achieves more
precise estimators.
\end{exm}

\begin{exm}[(A simple model)]\label{ex3}
To illustrate the adaptivity of our
algorithm to high dimension, we consider the following simple single-index
model:
%
\begin{equation}\label{simulation-SSIM}
Y=({\bolds{\beta}}^{\top}{\mathbf{X}})^2+\varepsilon.
\end{equation}
The true parameter is ${\bolds{\beta}}=(2,1,0,\ldots,0)^{\top}/\sqrt{5}
$; ${\mathbf{X}}$ is generated from $N_d(2,{\mathbf{I}})$. Both homogeneous errors and
heterogeneous ones are considered.
In the former case, $\varepsilon\sim N(0,0.2^2)$ and in
the latter case, $\varepsilon=\exp(\sqrt{5}{\bolds{\beta}}^{\top}{\mathbf{X}}/14)\widetilde{\varepsilon}$ with
$\widetilde{\varepsilon}\sim N(0,1)$. The latter case is designed
to show whether our method can handle heteroscedasticity.
A similar modeling setup was also used in \citet{WangXia2008}, Example 5.
The simulated results given in Table
\ref{simulation3} are based on 250
replicates with a sample of $n=100$ observations. An important
observation from this simulation is that the
proposed EFM approach still works even when the dimension of the
parameter is equal to or slightly larger than the number of
\begin{table}
\caption{Average estimation errors $\sum_{s=1}^d
|\hat{\beta}_s-\beta_s|$ for model (\protect\ref{simulation-SSIM})}\label{simulation3}
\begin{tabular*}{\tablewidth}{@{\extracolsep{\fill}}lccccc@{}}
\hline
$\bolds{\varepsilon}$ & & $\bolds{d=10}$ & $\bolds{d=50}$ & $\bolds{d=100}$ & $\bolds{d=120}$\\
\hline
& rMAVE & 0.0318 & 0.3484 & --- & --- \\
$\varepsilon\sim N(0,0.2^2)$& EDR & 0.0363 & 0.5020 & --- & --- \\
& EFM & 0.0272 & 0.2302 & 2.9409 & 5.0010 \\[5pt]
& rMAVE & 0.3427 & 4.6190 & --- & --- \\
$\varepsilon\sim N(0, \exp(\frac{2X_1+X_2}{7}))$& EDR & 0.2542 & 2.1112 & --- & ---
\\
& EFM & 0.2201 & 1.7937 & 4.1435 & 6.4973\\
\hline
\end{tabular*}
\legend{--- means that the
values cannot be calculated by rMAVE and EDR because of high dimension.}
\end{table}%
observations. It can be seen from Table
\ref{simulation3} that our approach also performs well under the
heteroscedasticity setup.
\end{exm}

\begin{exm}[(An oscillating function model)]\label{ex4}
A single-index model is designed as
%
\begin{equation}\label{simulation-OSIM}
Y=\sin(a{\bolds{\beta}}^{\top}{\mathbf{X}})+\varepsilon,
\end{equation}
where ${\bolds{\beta}}=(2,1,0,\ldots,0)^{\top}/\sqrt{5} $, ${\mathbf{X}}$ is
generated
from $N_d(2,{\mathbf{I}})$ and $\varepsilon\sim N(0,0.2^2)$. The number
of replications is 250 and the sample size $n=400$. The simulation
results are shown in Table~\ref{simulation4}. In these chosen values
for $a$, we see that EFM performs better than rMAVE and EDR. But as
\begin{table}[b]
\caption{Average estimation errors $\sum_{s=1}^d
|\hat{\beta}_s-\beta_s|$ for model (\protect\ref{simulation-OSIM})}\label{simulation4}
\begin{tabular}{@{}lcccccc@{}} 
\hline
& \multicolumn{3}{c}{$\bolds{a=\pi/2}$} &\multicolumn{3}{c@{}}{$\bolds{a=3\pi/4}$}
\\[-5pt]
&\multicolumn{3}{c}{\hrulefill}&\multicolumn{3}{c@{}}{\hrulefill}\\
$\bolds{d}$ & \textbf{rMAVE} & \textbf{EDR} & \textbf{EFM}& \textbf{rMAVE} & \textbf{EDR} &
\textbf{EFM}\\
\hline
10 & 0.0981 & 0.0918 & 0.0737 & 0.0970 & 0.0745 & 0.0725 \\
50 & 0.5247 & 0.6934 & 0.4355 & 0.6350 & 1.8484 & 0.5407 \\
\hline
\end{tabular}
\end{table}%
is understood, more oscillating functions are more difficult to handle
than those less oscillating functions.\vadjust{\goodbreak}
\end{exm}

\begin{exm}[(Comparison of variance)]\label{ex5}
To make our simulation results comparable with those of \citet{Carroll1997},
we mimic their simulation setup. Data of size $200$ are generated
according to the following model:
%
\begin{equation}\label{simulation-var}
Y_i=\sin\{\pi({\bolds{\beta}}^{\top}{\mathbf{X}}_i-A)/(B-A)\} + \alpha Z_i +
\varepsilon_i,
\end{equation}
where ${\mathbf{X}}_i$ are trivariate with independent $U(0,1)$
components, $Z_i$ are independent of ${\mathbf{X}}_i$ and $Z_i=0$ are for
$i$ odd and $Z_i=1$ for $i$ even, and
$\varepsilon_i$ follow a normal distribution $N(0,0.01)$ independent
of both ${\mathbf{X}}_i$ and $Z_i$. The
parameters are taken to be ${\bolds{\beta}}=(1,1,1)^{\top}/\sqrt{3}$,
$\alpha=0.3$, $A=\sqrt{3}/2-1.645/\sqrt{12}$ and
$B=\sqrt{3}/2+1.645/\sqrt{12}$.
Note that the EFM approach can still be applicable
for this model as the conditionally centered response $Y$
given $Z$ has the model as, because of the independence
between $\mathbf{X}$ and $Z$,
\[
Y_i-E(Y_i|Z_i)=a+ \sin\{\pi({\bolds{\beta}}^{\top}{\mathbf{X}}_i-A)/(B-A)\} +
\varepsilon_i.
\]
As $Z_i$ are dummy variables, estimating $E(Y_i|Z_i)$
is simple. Thus, when we regard $Y_i-E(Y_i|Z_i)$ as response,
the model is still a single-index model.
Here the number of replications is 100.
The method derived from \citet{Carroll1997} is referred to be the
GLPSIM approach. The numerical results are reported in Table
\ref{simulation5}. It shows that compared with the GPLSIM estimates,
the EFM estimates have smaller bias and smaller (or equal)
variance. Also in this example both EFM and GPLSIM can provide
reasonably accurate estimates.
\end{exm}

\begin{table}
\caption{Estimation for $\bolds{\beta}$ of model
(\protect\ref{simulation-var}) based on two randomly chosen samples}\label{simulation5}
\begin{tabular*}{\tablewidth}{@{\extracolsep{\fill}}lcccccc@{}}
\hline
& \multicolumn{3}{c}{\textbf{One group of sample}} &\multicolumn{3}{c@{}}{\textbf{Another group of sample}}\\[-5pt]
&\multicolumn{3}{c}{\hrulefill}&\multicolumn{3}{c@{}}{\hrulefill}\\
& $\bolds{X_1}$ & $\bolds{X_2}$ & $\bolds{X_3}$ & $\bolds{X_1}$ & $\bolds{X_2}$ &
$\bolds{X_3}$\\
\hline
GPLSIM est.
&0.595\tabnoteref{table5} &0.568\tabnoteref{table5} &0.569\tabnoteref{table5}
&0.563\tabnoteref{table5} &0.574\tabnoteref{table5}
&0.595\tabnoteref{table5}\\
GPLSIM s.e. &0.013\tabnoteref{table5} &0.013\tabnoteref{table5} &0.013\tabnoteref{table5}
&0.010\tabnoteref{table5} &0.010\tabnoteref{table5}
&0.010\tabnoteref{table5}\\
EFM est.
&0.579\hphantom{\tabnoteref{table5}}&0.575\hphantom{\tabnoteref{table5}}&0.577\hphantom{\tabnoteref{table5}} &0.573\hphantom{\tabnoteref{table5}}&0.577\hphantom{\tabnoteref{table5}}&0.580\hphantom{\tabnoteref{table5}}\\
EFM s.e. &0.011\hphantom{\tabnoteref{table5}}&0.011\hphantom{\tabnoteref{table5}}&0.011\hphantom{\tabnoteref{table5}} &0.010\hphantom{\tabnoteref{table5}}&0.010\hphantom{\tabnoteref{table5}} &0.010\hphantom{\tabnoteref{table5}}\\
\hline
\end{tabular*}
\tabnotetext{table5}{The
values are adopted from \citet{Carroll1997}. We abbreviate
``estimator'' to ``est.'' and ``standard error'' to ``s.e.,''
which are computed from the sample version of
${\bolds{\Sigma}}_{\hat{\bolds{\beta}}}$
defined in (\ref{asym-var}).}
\end{table}

\textit{Performance of profile quasi-likelihood ratio test}. To
illustrate how the profile quasi-likelihood ratio performs for linear
hypothesis problems, we simulate the same data as above, except that we allow
some components of the index to follow the null hypothesis:
\[
H_0\dvtx \beta_4=\beta_5=\cdots=\beta_d=0.
\]
We examine the power of the test under a sequence of the alternative hypotheses
indexed by parameter $\delta$ as follows:
\[
H_1\dvtx  \beta_4=\delta, \qquad \beta_s=0 \qquad \mbox{for } s\geq5.
\]
When $\delta=0$, the alternative hypothesis becomes the null hypothesis.

\begin{figure}

\includegraphics{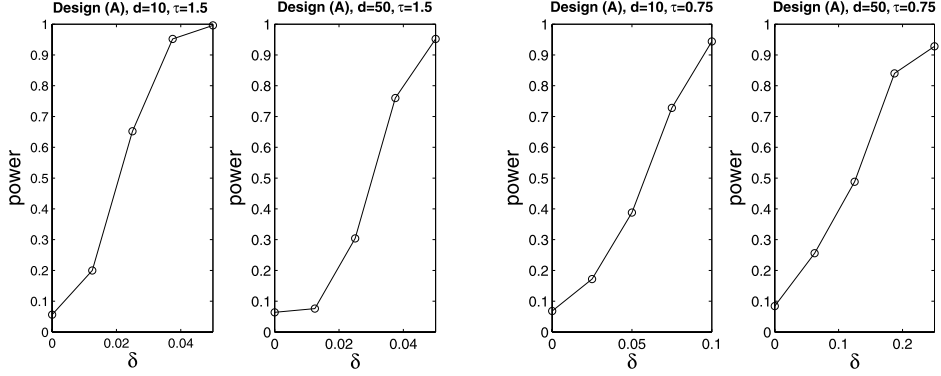}

  \caption{Simulation results for Design \textup{(A)} in Example \protect\ref{ex1}. The
left graphs depict the case $\tau=1.5$ with $\tau$ the first
parameter in $\operatorname{Beta}(\tau,1)$. The right graphs are for
$\tau=0.75$.}\label{fig1}
\end{figure}

\begin{figure}[b]

\includegraphics{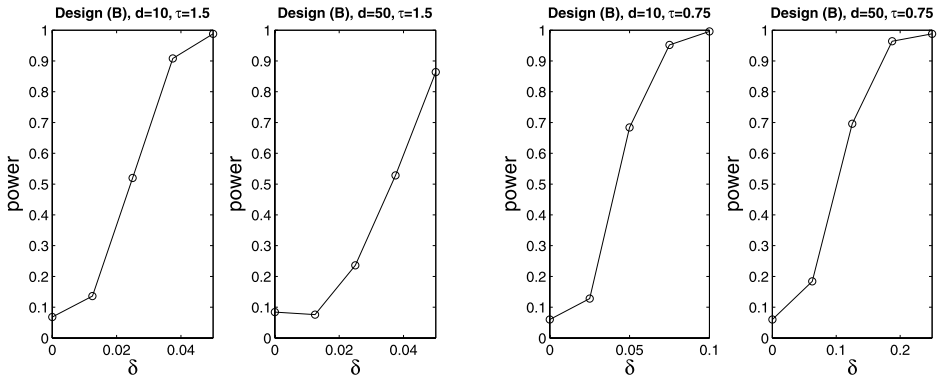}

  \caption{Simulation results for Design \textup{(B)} in Example \protect\ref{ex1}. The
left graphs depict the case $\tau=1.5$ with $\tau$ the first
parameter in $\operatorname{Beta}(\tau,1)$. The right graphs are for $\tau=0.75$.}\label{fig2}
\end{figure}

\begin{figure}

\includegraphics{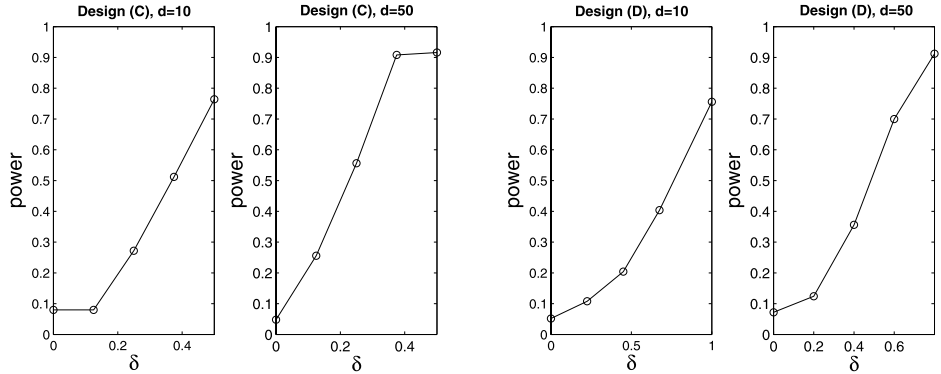}

  \caption{Simulation results for Example \protect\ref{ex2}. The left graphs
depict the case of Design \textup{(C)} with parameter dimension being 10 and
50. The right graphs are for Design \textup{(D)}.}\label{fig3}
\end{figure}

We examine the profile quasi-likelihood ratio test under a sequence
of alternative models, progressively deviating from the null
hypothesis, namely, as $\delta$ increases. The power functions are
calculated at the significance level: $0.05$, using the asymptotic
distribution. We calculate test statistics from 250 simulations by
employing the fixed-point algorithm and find the percentage of test
statistics greater than or equal to the associated quantile of the
asymptotic distribution. The pictures in Figures~\ref{fig1},~\ref{fig2} and~\ref{fig3}
illustrate the power function curves for two models under the given
significance levels. The power curves increase rapidly with~$\delta$, which shows the profile quasi-likelihood ratio test is
powerful. When $\delta$ is close to~0, the test sizes are all
approximately the significance levels.

\subsection{A real data example}\label{sec3.3}

Income, to some extent, is considered as an index of a successful life.
It is generally believed that demographic information, such as education level,
relationship in the household, marital status, the fertility rate and gender,
among others, has effects on amounts of income. For example,
\citet{Murray1997} illustrated that adults with higher intelligence have
higher income. \citet{ICKDDM1996} predicted income using a Bayesian
classifier offered by a machine learning algorithm. \citet{Madalozzo2008}
examined income differentials between married women and those who
remain single
or cohabit by using multivariate linear regression. Here we will use the
single-index model to explore the relationship between income and
some of its possible determinants.

We use the ``Adult'' database, which was extracted from the
Census Bureau database and is available on website:
\texttt{
\href{http://archive.ics.uci.edu/ml/datasets/Adult}{http://archive.ics.uci.edu/}
\href{http://archive.ics.uci.edu/ml/datasets/Adult}{ml/datasets/Adult}}. It was originally used to
model income exceeds over
USD 50,000$/$year based on census data. The purpose of using this
example is to understand the personal income patterns and demonstrate
the performance of the EFM method in real data analysis.
After excluding a few missing
data, the data set in our study includes 30,162 subjects. The
selected explanatory variables are:
\begin{itemize}
\item \textit{sex} (categorical): 1${}={}$Male, 0${}={}$Female.
\item \textit{native-country} (categorical): 1${}={}$United-States, 0${}={}$others.
\item \textit{work-class} (categorical): 1${}={}$Federal-gov, 2${}={}$Local-gov, 3${}={}$Private, 4${}={}$Self-emp-inc (self-employed, incorporated),
5${}={}$Self-emp-not-inc (self-\break employed, not incorporated), 6${}={}$State-gov.
\item \textit{marital-status} (categorical): 1${}={}$Divorced,
2${}={}$Married-AF-spouse (married, armed forces spouse present),
3${}={}$Married-civ-spouse (married, civilian spouse present),
4${}={}$Married-spouse-absent [married, spouse absent (exc.
separated)], 5${}={}$Never-married, 6${}={}$Separated, 7${}={}$Widowed.
\item \textit{occupation} (categorical): 1${}={}$Adm-clerical
(administrative support and clerical), 2${}={}$Armed-Forces, 3${}={}$Craft-repair, 4${}={}$Exec-managerial (execu\-tive-managerial), 5${}={}$Farming-fishing, 6${}={}$Handlers-cleaners, 7${}={}$Machine-op-inspct (machine
operator inspection), 8${}={}$Other-service, 9${}={}$Priv-house-serv (private
household services), 10${}={}$Prof-specialty (professional specialty), 11${}={}$Protective-serv, 12${}={}$Sales, 13${}={}$Tech-support, 14${}={}$Trans\-port-moving.
\item \textit{relationship} (categorical): 1${}={}$Husband, 2${}={}$Not-in-family, 3${}={}$Other-rela\-tive, 4${}={}$Own-child, 5${}={}$Unmarried, 6${}={}$Wife.
\item \textit{race} (categorical): 1${}={}$Amer-Indian-Eskimo, 2${}={}$Asian-Pac-Islander, 3${}={}$Black, 4${}={}$Other, 5${}={}$White.
\item \textit{age} (integer): number of years of age and greater
than or equal to 17.
\item \textit{fnlwgt} (continuous): The final sampling weights on
the CPS files are controlled to independent estimates of the civilian
noninstitutional population of the United States.
\item \textit{education} (ordinal): 1${}={}$Preschool (less than 1st
Grade), 2${}={}$1st--4th, 3${}={}$5th--6th, 4${}={}$7th--8th, 5${}={}$9th, 6${}={}$10th, 7${}={}$11th, 8${}={}$12th (12th Grade no Diploma), 9${}={}$HS-grad (high school
Grad-Diploma or Equiv), 10${}={}$Some-college (some college but no degree),
11${}={}$Assoc-voc (associate degree-occupational/vocational), 12${}={}$Assoc-acdm (associate degree-academic\break \mbox{program}),
13${}={}$Bachelors, 14${}={}$Masters, 15${}={}$Prof-school (professional school), 16${}={}$Doctorate.
\item \textit{education-num} (continuous): Number of years of education.
\item \textit{capital-gain} (continuous): A profit that results from
investments into a capital asset.
\item \textit{capital-loss} (continuous): A loss that results from
investments into a capital asset.
\item \textit{hours-per-week} (continuous): Usual number of hours
worked per week.
\end{itemize}

Note that all the explanatory variables up to ``age'' are categorical
with more than
two categories. As such, we use dummy variables to link up the
corresponding categories.
Specifically, for every original explanatory variable up to ``age,''
we use dummy variables to indicate it in which the number of dummy variables
is equal to the number of categories minus one.
By doing so, we then have 41 explanatory variables, where the first 35 ones
are dummy and the remaining ones are continuous.
After a preliminary data check, we find that the explanatory
variables $X_{37}=\mbox{``fnlwgt,''}$ $X_{39}=\mbox{``capital-gain''}$
and $X_{40}=\mbox{``capital-loss''}$ are very skewed to the left and
the latter two often take zero value. Before fitting
(\ref{real-data-model}) we first make a logarithm transformation for
these three variables
to have $\log(\mbox{``fnlwgt''})$,
$\log(1+\mbox{``capital-gain''})$ and
$\log(1+\mbox{``capital-loss''})$.
To make the explanatory variables
comparable in scale, we standardize each of them individually to obtain
mean 0
and variance 1. Since ``education'' and ``education-num''
are correlated, ``education'' is dropped from the model and
it results in a significantly smaller mean residual deviance.

\begin{table}
\tablewidth=320pt
\caption{Fitted coefficients for model (\protect\ref{real-data-model})
(estimated standard errors in parentheses)}\label{real-data}
\begin{tabular*}{\tablewidth}{@{\extracolsep{\fill}}ld{2.12}d{2.12}@{}}
\hline
\textbf{Variables} & \multicolumn{1}{c}{$\hat{\bolds{\beta}}$ \textbf{of SIM}} & \multicolumn{1}{c@{}}{$\hat{\bolds{\beta}}$ \textbf{of LR}}\\
\hline
Sex & 0.1102 \ (0.0028) & 0.1975 \ (0.0181)\\[3pt]
Native-country & 0.0412 \ (0.0027) &0.0354 \ (0.0116) \\[3pt]
Work-class\\
\quad Federal-gov &0.1237 \ (0.0059) &0.0739 \ (0.0108) \\
\quad Local-gov &0.2044 \ (0.0065) &0.0155 \ (0.0135) \\
\quad Private &-0.2603 \ (0.0075) &0.0775 \ (0.0200) \\
\quad Self-em-inc &0.1252 \ (0.0068) &0.0520 \ (0.0112) \\
\quad Self-emp-not-inc &0.1449 \ (0.0066) &-0.0157 \ (0.0147) \\[3pt]
Marital-Status\\
\quad Divorced &-0.0353 \ (0.0061) &-0.0304 \ (0.0264) \\
\quad Married-AF-spouse &0.0195 \ (0.0036) &0.0333 \ (0.0079) \\
\quad Married-civ-spouse &0.3257 \ (0.0150) &0.4545 \ (0.0754) \\
\quad Married-spouse-absent &-0.0115 \ (0.0029) &-0.0095 \ (0.0146) \\
\quad Never-married & -0.1876 \ (0.0085)&-0.1452 \ (0.0370) \\
\quad Separated &-0.0412 \ (0.0050) &-0.0221 \ (0.0179) \\[3pt]
Occupation\\
\quad Adm-clerical &-0.0302 \ (0.0050) &0.0131 \ (0.0164) \\
\quad Armed-Forces &-0.0086 \ (0.0031) &-0.0091 \ (0.0131) \\
\quad Craft-repair &-0.0913 \ (0.0050) &0.0263 \ (0.0146) \\
\quad Exec-managerial &0.1813 \ (0.0061) &0.1554 \ (0.0148) \\
\quad Farming-fishing &-0.0370 \ (0.0036) &-0.0772 \ (0.0125) \\
\quad Handlers-cleaners &-0.0947 \ (0.0033) &-0.0662 \ (0.0153) \\
\quad Machine-op-inspct &-0.1067 \ (0.0038) &-0.0290 \ (0.0133) \\
\quad Other-service &-0.1227 \ (0.0045) &-0.1192 \ (0.0195) \\
\quad Priv-house-serv &-0.0501 \ (0.0020) &-0.0833 \ (0.0379) \\
\quad Prof-specialty &0.2502 \ (0.0065) &0.1153 \ (0.0160) \\
\quad Protective-serv &0.1954 \ (0.0061) &0.0508 \ (0.0095) \\
\quad Sales &0.0316 \ (0.0050) &0.0615 \ (0.0147) \\
\quad Tech-support &0.0181 \ (0.0037)&0.0619 \ (0.0102) \\[3pt]
Relationship\\
\quad Husband &-0.1249 \ (0.0093) &-0.3264 \ (0.0254) \\
\quad Not-in-family &-0.0932 \ (0.0093) &-0.2074 \ (0.0612) \\
\quad Other-relative &-0.0958 \ (0.0038) &-0.1498 \ (0.0219) \\
\quad Own-child &-0.2218 \ (0.0076) &-0.3769 \ (0.0498) \\
\quad Unmarried &-0.1124 \ (0.0067) & -0.1739 \ (0.0446)\\[3pt]
Race\\
\quad Amer-Indian-Eskimo &-0.0252 \ (0.0024) & -0.0226
\ (0.0109)\\
\quad Asian-Pac-Islander &0.0114 \ (0.0030) &0.0062 \ (0.0101) \\
\quad Black &-0.0300 \ (0.0024) &-0.0182 \ (0.0111) \\
\quad Other &-0.0335 \ (0.0021) &-0.0286 \ (0.0129) \\
\hline
\end{tabular*}
\end{table}
\setcounter{table}{5}
\begin{table}
\tablewidth=279pt
\caption{(Continued)}
\begin{tabular*}{\tablewidth}{@{\extracolsep{\fill}}ld{2.12}d{2.13}@{}}
\hline
\textbf{Variables} & \multicolumn{1}{c}{$\hat{\bolds{\beta}}$ \textbf{of SIM}} & \multicolumn{1}{c@{}}{$\hat{\bolds{\beta}}$ \textbf{of LR}}\\
\hline
Age &0.2272 \ (0.0042) &0.1798 \ (0.0111) \\
Fnlwgt &0.0099 \ (0.0028) &0.0414 \ (0.0092) \\
Education-num&0.4485 \ (0.0045) &0.3732 \ (0.0122) \\
Capital-gain &0.2859 \ (0.0055) &0.2582 \ (0.0084) \\
Capital-loss &0.1401 \ (0.0042) &0.1210 \ (0.0078) \\
Hours-per-week &0.2097 \ (0.0035) &0.1823 \ (0.0101) \\
\hline
\end{tabular*}
\end{table}

The single-index model will be used to model the relationship
between income and the relevant 43 predictors ${\mathbf{X}}=(X_1,\ldots,X_{43})^{\top} $:
%
\begin{equation}\label{real-data-model}
P(\mbox{``income''}>50\mbox{,}000|{\mathbf{X}})= \exp\{g({\bolds{\beta}}^{\top}{\mathbf{X}})\}/[1+\exp\{g({\bolds{\beta}}^{\top}{\mathbf{X}})\}],
\end{equation}
where $Y=I(\mbox{``income''}>50\mbox{,}000)$ and ${\bolds{\beta}}=(\beta
_1,\ldots,\beta_{43})^{\top}$ and $\beta_s$
represents the effect of the $s$th predictor. Formally, we are
testing the effect of gender, that is,
%
\begin{equation}\label{real1}
H_0\dvtx\beta_1=0 \quad \longleftrightarrow\quad  H_1\dvtx\beta_1\neq0.
\end{equation}

The fixed-point iterative
algorithm is employed to compute the estimate for $\bolds{\beta}$. To
illustrate further the practical implications of this approach, we
compare our results to those obtained by using an ordinary logistic
regression (LR). The coefficients of the two models are given in Table
\ref{real-data}. To make the analyses presented in the table
comparable, we consider two standardizations. First, we standardize
every explanatory variable with mean 0 and variance 1 so that
the coefficients can be used to compare the relative influence
from different explanatory variables. However, such a standardization
does not allow us to compare between the single-index model and the
ordinary logistic regression model. We then further
normalize the coefficients to be with Euclidean norm 1, and then the
estimates of their standard
errors are also adjusted accordingly.
The single-index model provides more
reasonable results: $X_{38}=\mbox{``education-num''}$ has its strongest
positive effect on income; those who got a bachelor's degree or
higher seem to have much higher income than those with lower
education level. In contrast, results derived from a logistic
regression show
that ``married-civ-spouse'' is the largest positive contributor.

%
%

Some other interesting conclusions could be obtained by looking at
the output. Both ``sex'' and ``native-country'' have a positive effect.
Persons
who worked without pay in a family business,
unpaid childcare and others earn a lower income than persons
who worked for wages or for themselves. The
``fnlwgt'' attribute has a positive relation to income.
Males are likely to
make much more money than females.
The expected sign for
marital status except the {\it married} (married-AF-spouse,
married-civ-spouse) is negative, given that the household production
theory affirms that division of work is efficient when each member of a family
dedicates his or her time to the more productive job. Men usually
receive relatively better compensation for their time in the
labor market than in home production. Thus, the expectation is
that married women dedicate more time to home tasks and less to
the labor market, and this would imply a different probability
of working given the marital status choice.

Also ``race'' influences the income and Asian
or Pacific Islanders seem to make more money than other races.
And also, one's income significantly increases as working hours
increase. Both ``capital-gain'' and ``capital-loss'' have
positive effects, so we think that people make more money who
can use more money to invest. The presence of young children
has a negative influence on the income. ``age'' accounts for the
experience effect and has a positive effect. Hence the conclusion based
on the single-index model is consistent with what we expect.

\begin{figure}

\includegraphics{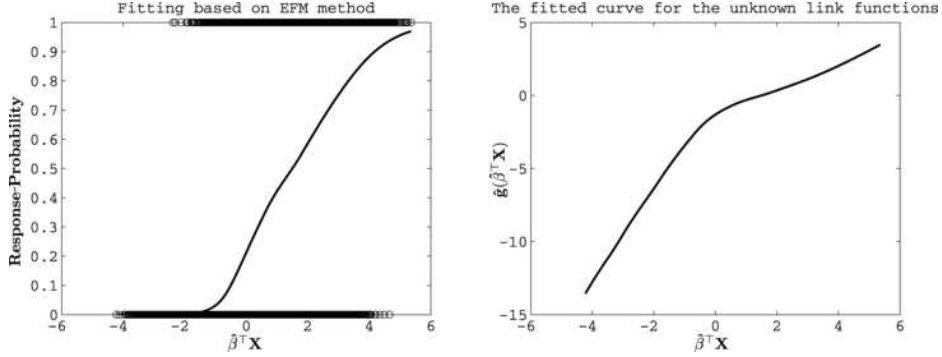}

  \caption{Adult data: The left graph is a plot of
predicted response probability based on the single-index
model. The right graph is the fitted curve for the unknown link
function $g(\cdot)$.}\label{fig4}
\end{figure}

To help with
interpretation of the model, plots of ${\bolds{\beta}}^{\top}{\mathbf{X}}$
versus predicted response probability and
$\hat{g}({\bolds{\beta}}^{\top}{\mathbf{X}})$ are generated, respectively,
and can be found on the right column in Figure~\ref{fig4}.
When the estimated single-index is greater than 0,
$\hat{g}(\hat{\bolds{\beta}}{\mathbf{X}})$ shows some degree of curvature.
An alternative choice is to fit the data using generalized partially
linear additive models (GPLAM) with nonparametric components of continuous
explanatory variables. The relationships among ``age,'' ``fnlwgt,''
``capital-gain,''
``capital-loss'' and ``hours-per-week'' all show nonlinearity. The
mean residual deviances of SIM, LR and GPLAM are $0.7811$, $0.6747$ and
$0.6240$, respectively. SIM under study
provides a slightly worse fit than the others.
However, we note that LR is, up to a link function, linear about $\mathbf{X}$, and,
according to the results of GPLAM, which is a more general model than
LR, the actual
relationship cannot have such a structure. SIM can reveal nonlinear structure.
On the other hand, although the minimum mean residual deviance
can be not surprisingly attained by GPLAM, this model has,
respectively,
$\approx$ 34 and 41 more degrees of freedom
than SIM and LR have.

We now employ the quasi-likelihood ratio test to the test problem
(\ref{real1}). The QLR test statistic is $166.52$ with one degree
of freedom, resulting in a $P$-value of $<10^{-5}$. Hence this result
provides strong evidence that gender has a significant influence
on high income.

The Adult data set used in this paper is a rich data set.
Existing work mainly focused on the prediction accuracy based
on machine learning methods. We make an attempt to explore the
semiparametric regression pattern suitable for the data. Model
specification and
variable selection merit further study.

\begin{appendix}

\section*{Appendix: Outline of proofs}\label{app}
\setcounter{equation}{0}

We first introduce some regularity conditions.

\textit{Regularity Conditions}:
\begin{enumerate}[(a)]
\item[(a)]$\mu(\cdot), V(\cdot), g(\cdot), {\mathbf{h}}(\cdot)=E({\mathbf{X}}|
{\bolds{\beta}}^{\top} {\mathbf{X}}=\cdot)$ have two bounded and continuous
derivatives. $V(\cdot)$ is uniformly bounded and
bounded away from 0.
\item[(b)] Let
$q(z,y)=\mu'(z)V^{-1}(z)\{y-\mu(z)\}$. Assume that $\partial
q(z,y)/\partial z<0$ for $z\in\mathbb{R}$ and $y$ in the range of
the response variable.
\item[(c)] The largest eigenvalue of ${\bolds{\Omega}}_{22}$ is bounded away
from infinity.
\item[(d)] The density function $f_{{\bolds{\beta}}^{\top}{\mathbf{x}}}({\bolds{\beta}}^{\top} {\mathbf{x}})$ of random variable
${\bolds{\beta}}^{\top}{\mathbf{X}}$ is bounded away from 0 on $T_{\bolds{\beta}}$
and satisfies the Lipschitz condition of order 1 on $T_{\bolds{\beta}}$,
where $T_{\bolds{\beta}}=\{{\bolds{\beta}}^{\top} {\mathbf{x}}\dvtx  {\mathbf{x}}\in T\}$ and
$T$ is a compact support set of ${\mathbf{X}}$.
\item[(e)] Let $Q^*[{\bolds{\beta}}]=
\int Q[\mu\{g({\bolds{\beta}}^{\top}{\mathbf{x}})\},y]
f(y|{\bolds{\beta}}^{0\top}{\mathbf{x}})f({\bolds{\beta}}^{0\top}{\mathbf{x}})\, dy\,d({\bolds{\beta}}^{0\top}{\mathbf{x}})$ with ${\bolds{\beta}}^0$ denoting
the true parameter value and
$Q[\mu,y]=\int_{\mu}^y\frac{s-y}{V\{\mu^{-1}(s)\}}\,ds$.\vspace*{-2pt} Assume that
$Q^*[{\bolds{\beta}}]$ has a unique maximum at
${\bolds{\beta}}={\bolds{\beta}}^0$, and
\[
E \Bigl[
\sup_{{\bolds{\beta}}^{(1)}}\sup_{{\bolds{\beta}}^{\top}{\mathbf{X}}} |\mu'\{g({\bolds{\beta}}^{\top}{\mathbf{X}})\}
V^{-1}\{g({\bolds{\beta}}^{\top}{\mathbf{X}})\}[Y-\mu\{g({\bolds{\beta}}^{\top}{\mathbf{X}})\}] |^2 \Bigr] < \infty
\]
and $E\|{\mathbf{X}}\|^2< \infty$.\vadjust{\goodbreak}
\item[(f)] The kernel $K$ is a bounded and symmetric density function with a bounded
derivative, and satisfies
\[
\int_{-\infty}^{\infty} t^2K(t)\,dt\neq0 \quad \mbox{and}\quad
\int_{-\infty}^{\infty} |t|^jK(t)\,dt<\infty,\qquad  j=1,2,\ldots.
\]
\end{enumerate}

Condition (a) is some mild smoothness conditions on the involved
functions of the model. We impose condition (b) to guarantee that
the solutions of (\ref{eq-kernel}), $\hat{g}(t)$ and
$\hat{g}'(t)$, lie in a compact set. Condition (c) implies that the
second moment of estimating equation (\ref{peqp}), $\operatorname{tr}({\mathbf{J}}^{\top} {\bolds{\Omega}} {\mathbf{J}})$, is bounded. Then the CLT can be
applied to $G({\bolds{\beta}})$. Condition (d) means that ${\mathbf{X}}$ may
have discrete components and the density function of
${\bolds{\beta}}^{\top}{\mathbf{X}}$ is positive, which ensures that the
denominators involved in the nonparametric estimators, with high
probability, are bounded away from 0. The
uniqueness condition in condition (e) can be checked in the following
case for example. Assume that $Y$ is a Poisson variable with mean
$\mu\{g({\bolds{\beta}}^{\top}{\mathbf{x}})\} = \exp\{g({\bolds{\beta}}^{\top
}{\mathbf{x}})\}$.
The maximizer $\beta_0$ of $ Q^*[{\bolds{\beta}}]$ is equal to\vspace*{1pt}
the solution of the equation $E [E \{[\exp\{g({\bolds{\beta}}^{0\top}
{\mathbf{X}})\}-\exp\{g({\bolds{\beta}}^{\top}
{\mathbf{X}})\}]g'({\bolds{\beta}}^{\top}
{\mathbf{X}})\}{\mathbf{J}}^{\top}{\mathbf{X}} |{\bolds{\beta}}^{0\top}{\mathbf{X}} \} ]=0$.
${\bolds{\beta}}_0$\vspace*{1pt} is unique
when $g'(\cdot)$ is not a zero-valued constant function and
the matrix ${\mathbf{J}}^{\top}E({\mathbf{X}}{\mathbf{X}}^{\top}){\mathbf{J}}$ is not
singular. Under the second part
of condition (e),
it is permissible to interchange differentiation and integration
when differentiating $E [Q[\mu\{g({\bolds{\beta}}^{\top} {\mathbf{X}})\},Y] ]$. Condition (f) is a commonly used smoothness
condition, including the Gaussian kernel and the quadratic kernel.
All of the conditions can be relaxed at the expense of longer
proofs.

Throughout the \hyperref[app]{Appendix}, $Z_n=\CO_P(a_n)$ denotes that $a_n^{-1}Z_n$
is bounded
in probability and the derivation for the order of $Z_n$ is based on
the fact
that $Z_n=\CO_P\{\sqrt{E(Z_n^2)}\}$. Therefore, it allows to apply the
Cauchy--Schwarz inequality to the quantity having stochastic order
$a_n$.\vspace*{-2pt}

\subsection{\texorpdfstring{Proof of Proposition \protect\ref{prop-non}}{Proof of Proposition 1}}\label{sec4.1}
We outline the proof here, while the details are given in the
supplementary materials
\citep{CuiHardleZhu2010s}.

\begin{longlist}[(iii)]
\item[(i)] Conditions (a), (b), (d) and (f) are essentially equivalent
conditions given by \citet{Carroll1998}, and as a consequence the
derivation of bias and variance for $\hat{g}({\bolds{\beta}}^{\top} {\mathbf{x}})$ and $\hat{g}'({\bolds{\beta}}^{\top}{\mathbf{x}})$ is similar to that of
\citet{Carroll1998}.

\item[(ii)] The first equation of (\ref{eq-kernel}) is
\begin{eqnarray*}
0&=&\sum_{j=1}^n K_h({\bolds{\beta}}^{\top}
{\mathbf{X}}_j-{\bolds{\beta}}^{\top}
{\mathbf{x}}) \mu'\{\hat{\alpha}_0+\hat{\alpha}_1({\bolds{\beta}}^{\top}{\mathbf{X}}_j-{\bolds{\beta}}^{\top} {\mathbf{x}})\}\\[-3pt]
&&\hphantom{\sum_{j=1}^n}{} \times V^{-1}\{\hat{\alpha}_0+\hat{\alpha}_1({\bolds{\beta}}^{\top
}{\mathbf{X}}_j-{\bolds{\beta}}^{\top}
{\mathbf{x}})\}[Y_j-\mu\{\hat{\alpha}_0+\hat{\alpha}_1({\bolds{\beta}}^{\top
}{\mathbf{X}}_j-{\bolds{\beta}}^{\top} {\mathbf{x}})\}].
\end{eqnarray*}
Taking derivatives with respect to ${\bolds{\beta}}^{(1)}$ on both sides, direct
observations lead to
\[
\frac{\partial\hat{\alpha}_0}{\partial{\bolds{\beta}}^{(1)}} =
\{B({\bolds{\beta}}^{\top}{\mathbf{x}})\}^{-1}\{A_1({\bolds{\beta}}^{\top}
{\mathbf{x}})+A_2({\bolds{\beta}}^{\top}
{\mathbf{x}})+A_3({\bolds{\beta}}^{\top}{\mathbf{x}})\},\vadjust{\goodbreak}
\]
where
\begin{eqnarray*}
B({\bolds{\beta}}^{\top} {\mathbf{x}}) &=& -\sum_{j=1}^n
K_h({\bolds{\beta}}^{\top}{\mathbf{X}}_j-{\bolds{\beta}}^{\top}
{\mathbf{x}})q'_z\{\hat{\alpha}_0+\hat{\alpha}_1({\bolds{\beta}}^{\top}
{\mathbf{X}}_j-{\bolds{\beta}}^{\top}
{\mathbf{x}}),Y_j\},\\
A_1({\bolds{\beta}}^{\top}{\mathbf{x}}) &=& \sum_{j=1}^n
K_h({\bolds{\beta}}^{\top}{\mathbf{X}}_j-{\bolds{\beta}}^{\top} {\mathbf{x}}){\mathbf{J}}^{\top}
({\mathbf{X}}_j-{\mathbf{x}})q'_z\{\hat{\alpha}_0+\hat{\alpha}_1({\bolds{\beta}}^{\top}{\mathbf{X}}_j-{\bolds{\beta}}^{\top}
{\mathbf{x}}),Y_j\}\hat{\alpha}_1, \\
A_2({\bolds{\beta}}^{\top}{\mathbf{x}})&=& \sum_{j=1}^n
K_h({\bolds{\beta}}^{\top}{\mathbf{X}}_j-{\bolds{\beta}}^{\top} {\mathbf{x}})q'_z\{\hat{\alpha}_0+\hat{\alpha}_1({\bolds{\beta}}^{\top} {\mathbf{X}}_j-{\bolds{\beta}}^{\top}{\mathbf{x}}),Y_j\}\\
&&\hphantom{\sum_{j=1}^n}{}\times({\bolds{\beta}}^{\top} {\mathbf{X}}_j-{\bolds{\beta}}^{\top}{\mathbf{x}})\frac{\partial\hat{\alpha}_1}{
\partial{\bolds{\beta}}^{(1)}},\\
A_3({\bolds{\beta}}^{\top}
{\mathbf{x}}) &=& \sum_{j=1}^n h^{-1}
K'_h({\bolds{\beta}}^{\top}{\mathbf{X}}_j-{\bolds{\beta}}^{\top}{\mathbf{x}}){\mathbf{J}}^{\top}
({\mathbf{X}}_j-{\mathbf{x}})q\{\hat{\alpha}_0+\hat{\alpha}_1({\bolds{\beta}}^{\top
}{\mathbf{X}}_j-{\bolds{\beta}}^{\top} {\mathbf{x}}),Y_j\}
\end{eqnarray*}
with $K'_h(\cdot)=h^{-1}K'(\cdot/h)$. Note that $\partial
\hat{\alpha}_0/\partial{\bolds{\beta}}^{(1)}=\partial
\hat{g}({\bolds{\beta}}^{\top}{\mathbf{x}})/\partial{\bolds{\beta}}^{(1)}$; then
we have
%
\begin{eqnarray}\label{der-decom}
\frac{\partial\hat{g}({\bolds{\beta}}^{\top}
{\mathbf{x}})}{\partial{\bolds{\beta}}^{(1)}} &=&
\{B({\bolds{\beta}}^{\top}
{\mathbf{x}})\}^{-1}A_1({\bolds{\beta}}^{\top}{\mathbf{x}})\nonumber\\[-8pt]\\[-8pt]
&&{} +\{B({\bolds{\beta}}^{\top}{\mathbf{x}})\}^{-1}A_2({\bolds{\beta}}^{\top}{\mathbf{x}})+\{B({\bolds{\beta}}^{\top}{\mathbf{x}})\}^{-1}A_3({\bolds{\beta}}^{\top}{\mathbf{x}}).\nonumber
\end{eqnarray}
We will prove that
%
\begin{eqnarray}
&&E\|\{B({\bolds{\beta}}^{\top}{\mathbf{x}})\}^{-1}A_1({\bolds{\beta}}^{\top}{\mathbf{x}})
-g'({\bolds{\beta}}^{\top}{\mathbf{x}}){\mathbf{J}}^{\top}\{{\mathbf{x}}-{\mathbf{h}}({\bolds{\beta}}^{\top}
{\mathbf{x}})\}\|^2 \nonumber\\[-8pt]\\[-8pt]
&&\qquad = \CO_P(h^4+n^{-1}h^{-3}),\nonumber
\end{eqnarray}
the second term in (\ref{der-decom}) is of order $\CO
_P(h^4+n^{-1}h)$, and the
third term is of order $\CO_P(h^4+n^{-1}h^{-3})$. The combination of
(\ref{der-decom}) and these three results can directly lead to result
(ii) of
Proposition~\ref{prop-non}. The detailed proof is summarized in three
steps and is given in the supplementary materials \citep
{CuiHardleZhu2010s}.

\item[(iii)] By mimicking the proof of (ii), we can show
that (iii) holds. See supplementary materials for details.
\end{longlist}

\subsection{\texorpdfstring{Proofs of (\protect\ref{eq-profile}) and (\protect\ref{peqp})}{Proofs of (2.6) and (2.7)}}\label{sec4.2}
It is proved in the supplementary materials
\citep{CuiHardleZhu2010s}.

\subsection{\texorpdfstring{Proof of Theorem \protect\ref{asymptotic-normality}}{Proof of Theorem 2.1}}\label{sec4.3}

(i) Note that the estimating equation defined in
(\ref{eq-profile}) is just the gradient of the following
quasi-likelihood:
\[
\hat{Q}({\bolds{\beta}})=\sum_{i=1}^n
Q[\mu\{\hat{g}({\bolds{\beta}}^{\top} {\mathbf{X}}_i)\},Y_i]
\]
with $Q[\mu,y]=\int^{\mu} \frac{y-s}{V\{\mu^{-1}(s)\}}\,ds$ and
$\mu^{-1}(\cdot)$ is the inverse function of $\mu(\cdot)$. Then for
${\bolds{\beta}}^{(1)}$ satisfying $(\sqrt{1-\|{\bolds{\beta}}^{(1)}\|^2},
{\bolds{\beta}}^{(1)\top})^{\top}\in\Theta$, we have
\[
\hat{\bolds{\beta}}{}^{(1)}=\arg\max_{{\bolds{\beta}}^{(1)}}\hat{Q}({\bolds{\beta}}).
\]
The proof is based on Theorem 5.1 in \citet{Ichimura1993}. In that
theorem the consistency of ${\bolds{\beta}}^{(1)}$ is proved by means of
proving that
%
\begin{eqnarray}\label{uni1}
\sup_{{\bolds{\beta}}^{(1)}} \Biggl|\frac{1}{n}\sum_{i=1}^n
Q[\mu\{\hat{g}({\bolds{\beta}}^{\top} {\mathbf{X}}_i)\},Y_i]-\frac{1}{n}\sum_{i=1}^n
Q[\mu\{g({\bolds{\beta}}^{\top} {\mathbf{X}}_i)\},Y_i] \Biggr|&=&{
\mo}_P(1),
\\
\label{uni2}
\qquad \quad \sup_{{\bolds{\beta}}^{(1)}} \Biggl|\frac{1}{n}\sum_{i=1}^n
Q[\mu\{g({\bolds{\beta}}^{\top} {\mathbf{X}}_i)\},Y_i]-\frac{1}{n}\sum_{i=1}^n
E [Q[\mu\{g({\bolds{\beta}}^{\top} {\mathbf{X}}_i)\},Y_i] ] \Biggr|&=&{\mo}_P(1)
\end{eqnarray}
and
%
\begin{equation}\label{uni3}
\qquad \Biggl|\frac{1}{n}\sum_{i=1}^n
Q[\mu\{\hat{g}({\bolds{\beta}}_0^{\top} {\mathbf{X}}_i)\},Y_i]-\frac{1}{n}\sum_{i=1}^n
E [Q[\mu\{g({\bolds{\beta}}_0^{\top} {\mathbf{X}}_i)\},Y_i] ] \Biggr|={\mo}_P(1).
\end{equation}
Regarding the validity of (\ref{uni3}), this directly follows from
(\ref{uni1}) and (\ref{uni2}). The type of uniform convergence
result such as (\ref{uni2}) has been well established in the
literature; see, for example, \citet{Andrews1987}. We now verify the
validity of (\ref{uni1}), which reduces to showing the uniform
convergence of the estimator $\hat{g}(t)$ under condition (e)
[see \citet{Ichimura1993}]. This can be obtained in a similar way
as in \citet{Kong2010}, taking into account that the regularity
conditions imposed in Theorem~\ref{asymptotic-normality} are
stronger than the corresponding ones in that paper.

(ii) Recall the notation ${\mathbf{J}}, {\bolds{\Omega}}$ and ${\mathbf{G}}({\bolds{\beta}})$ introduced in Section~\ref{sec2}. By (\ref{peqp}), we have
shown that
%
\begin{equation}\label{asym-expan-beta}
\sqrt{n}\bigl(\hat{\bolds{\beta}}{}^{(1)}-{\bolds{\beta}}^{(1)0}\bigr)=\frac{1}{\sqrt
{n}} \{{\mathbf{J}}^{\top} {\bolds{\Omega}}{\mathbf{J}}\}^{+} {\mathbf{G}}({\bolds{\beta}}) +
{\mo}_P(1).
\end{equation}
Theorem~\ref{asymptotic-normality} follows directly from the above asymptotic
expansion and the fact that $E\{{\mathbf{G}}({\bolds{\beta}}){\mathbf{G}}^{\top}({\bolds{\beta}})\}=n{\mathbf{J}}^{\top}{\bolds{\Omega}}{\mathbf{J}}$.
\hfill$\square$

\subsection{\texorpdfstring{Proof of Corollary \protect\ref{variance-comparision}}{Proof of Corollary 1}}\label{sec4.4}

The asymptotic covariance of $\hat{\bolds{\beta}}$ can\vspace*{-1pt} be obtained by
adjusting the asymptotic covariance of $\hat{\bolds{\beta}}{}^{(1)}$ via
the multivariate delta method, and is of form ${\mathbf{J}}({\mathbf{J}}^{\top}{\bolds{\Omega}}{\mathbf{J}})^{+}{\mathbf{J}}^{\top}$. Next we will
compare this asymptotic covariance with that (denoted by
${\bolds{\Omega}}^+$) given in \citet{Carroll1997}. Write ${\bolds{\Omega}}$
as
\[
{\bolds{\Omega}}= \pmatrix{
{\bolds{\Omega}}_{11} & {\bolds{\Omega}}_{12} \cr
{\bolds{\Omega}}_{21} & {\bolds{\Omega}}_{22}
},
\]
where ${\bolds{\Omega}}_{22}$ is a $(d-1)\times(d-1)$ matrix. We will
next investigate two cases, respectively: $\det({\bolds{\Omega}}_{22})\neq
0$ and
$\det({\bolds{\Omega}}_{22})= 0$. Let ${\bolds{\alpha}}=
-{\bolds{\beta}}^{(1)}/\sqrt{1-\|{\bolds{\beta}}^{(1)}\|^2}=-{\bolds{\beta}
}^{(1)}/\beta_1$.

Consider the case that $\det({\bolds{\Omega}}_{22})\neq0$. Because
$\operatorname{rank}({\bolds{\Omega}})=d-1$, $\det({\bolds{\Omega}}_{11}{\bolds{\Omega}
}_{22}-{\bolds{\Omega}}_{21}{\bolds{\Omega}}_{12})=0$. Note that ${\bolds{\Omega}}_{22}$ is
nondegenerate; it can be easily shown that ${\bolds{\Omega}}_{11}={
\bolds{\Omega}}_{12}{\bolds{\Omega}}_{22}^{-1}{\bolds{\Omega}}_{21}$. Combining this
with the
following fact:
\begin{eqnarray*}
{\mathbf{J}}^{\top}{\bolds{\Omega}}{\mathbf{J}}&=& \pmatrix{
{\bolds{\alpha}} & {\mathbf{I}}_{d-1}
}\pmatrix{
{\bolds{\Omega}}_{11} & {\bolds{\Omega}}_{12} \cr
{\bolds{\Omega}}_{21} & {\bolds{\Omega}}_{22}}
\pmatrix{
{\bolds{\alpha}}^{\tau}\cr
{\mathbf{I}}_{d-1}
}\\
&=&{\bolds{\Omega}}_{22}+\bigl({\bolds{\Omega}}_{21}/\sqrt{{\bolds{\Omega}
}_{11}}+\sqrt{{\bolds{\Omega}}_{11}}{\bolds{\alpha}}\bigr)
\bigl({\bolds{\Omega}}_{12}/\sqrt{{\bolds{\Omega}}_{11}}+\sqrt{{\bolds{\Omega}
}_{11}}{\bolds{\alpha}}^{\top}\bigr)-{\bolds{\Omega}}_{21}
{\bolds{\Omega}}_{12}/{\bolds{\Omega}}_{11},
\end{eqnarray*}
we can get that ${\mathbf{J}}^{\top}{\bolds{\Omega}}{\mathbf{J}}$ is
nondegenerate. In this situation, its inverse\break $({\mathbf{J}}^{\top}{
\bolds{\Omega}}{\mathbf{J}})^{+}$ is just the ordinary inverse
$({\mathbf{J}}^{\top}{\bolds{\Omega}}{\mathbf{J}})^{-1}$. Then
${\mathbf{J}}({\mathbf{J}}^{\top}{\bolds{\Omega}}{\mathbf{J}})^{+}{\mathbf{J}}^{\top}=\break
\{{\mathbf{J}}({\mathbf{J}}^{\top}{\bolds{\Omega}}{\mathbf{J}})^{-1/2} \}
\{({\mathbf{J}}^{\top}\times{\bolds{\Omega}}{\mathbf{J}})^{-1/2}{\mathbf{J}}^{\top} \}$, a
full-rank decomposition. Then
\begin{eqnarray*}
\{{\mathbf{J}}({\mathbf{J}}^{\top}{\bolds{\Omega}}{\mathbf{J}})^{+}{\mathbf{J}}^{\top} \}^+
&= &
\{{\mathbf{J}}({\mathbf{J}}^{\top}{\bolds{\Omega}}{\mathbf{J}})^{-1/2} \}\\
&&{}\times\{({\mathbf{J}}^{\top}{\bolds{\Omega}}{\mathbf{J}})^{-1/2}{\mathbf{J}}^{\top
}{\mathbf{J}}({\mathbf{J}}^{\top}{\bolds{\Omega}}{\mathbf{J}})^{-1}{\mathbf{J}}^{\top}{\mathbf{J}}({\mathbf{J}}^{\top}{\bolds{\Omega}}{\mathbf{J}})^{-1/2} \}^{-1} \\
&&{}\times\{({\mathbf{J}}^{\top}{\bolds{\Omega}}{\mathbf{J}})^{-1/2}{\mathbf{J}}^{\top}
\}\\
&= &{\mathbf{J}}({\mathbf{J}}^{\top}{\mathbf{J}})^{-1}{\mathbf{J}}^{\top}{\bolds{\Omega}}{\mathbf{J}}({\mathbf{J}}^{\top}{\mathbf{J}})^{-1}{\mathbf{J}}^{\top}\\
&= & {\bolds{\Omega}}.
\end{eqnarray*}
This means that ${\mathbf{J}}({\mathbf{J}}^{\top}{\bolds{\Omega}}{\mathbf{J}})^{+}{\mathbf{J}}^{\top}
={\bolds{\Omega}}^{+}$.

When $\det({\bolds{\Omega}}_{22})= 0$, we can obtain that
\[
{\bolds{\Omega}}^{+}= \pmatrix{
1/{\bolds{\Omega}}_{11}+{\bolds{\Omega}}_{12}{\bolds{\Omega}}_{22.1}^+{
\bolds{\Omega}}_{21}/{\bolds{\Omega}}_{11}^2
& -
{\bolds{\Omega}}_{12}{\bolds{\Omega}}_{22.1}^+/{\bolds{\Omega}}_{11}\vspace*{3pt}\cr
- {\bolds{\Omega}}_{22.1}^+{\bolds{\Omega}}_{21}/{\bolds{\Omega}}_{11} &
{\bolds{\Omega}}_{22.1}^+ }
\]
with
${\bolds{\Omega}}_{22.1}={\bolds{\Omega}}_{22}-{\bolds{\Omega}}_{21}{\bolds{\Omega}
}_{12}/{\bolds{\Omega}}_{11}$.
Write ${\mathbf{J}}({\mathbf{J}}^{\top}{\bolds{\Omega}}{\mathbf{J}})^{+}{\mathbf{J}}^{\top}$ as
\[
\pmatrix{
{\bolds{\alpha}}^{\top}({\mathbf{J}}^{\top}{\bolds{\Omega}}{\mathbf{J}})^{+}{\bolds{\alpha}}
& {\bolds{\alpha}}^{\top}({\mathbf{J}}^{\top}{\bolds{\Omega}}{\mathbf{J}})^{+}\vspace*{2pt}
\cr
({\mathbf{J}}^{\top}{\bolds{\Omega}}{\mathbf{J}})^{+}{\bolds{\alpha}}
& ({\mathbf{J}}^{\top}{\bolds{\Omega}}{\mathbf{J}})^{+}
}.
\]
Note that ${\mathbf{J}}^{\top}{\bolds{\Omega}}{\mathbf{J}} =
{\bolds{\Omega}}_{22.1}+({\bolds{\Omega}}_{21}/\sqrt{{\bolds{\Omega}
}_{11}}+\sqrt{{\bolds{\Omega}}_{11}}{\bolds{\alpha}})
({\bolds{\Omega}}_{12}/\sqrt{{\bolds{\Omega}}_{11}}+\sqrt{{\bolds{\Omega}
}_{11}}{\bolds{\alpha}}^{\top})$,
so ${\mathbf{J}}^{\top}{\bolds{\Omega}}{\mathbf{J}} \geq{\bolds{\Omega}}_{22.1}.$
Combining this
with $\operatorname{rank}({\bolds{\Omega}}_{22})=d-2$, we have that $({\mathbf{J}}^{\top}{
\bolds{\Omega}}{\mathbf{J}})^{+} \leq{\bolds{\Omega}}_{22.1}^+$. It is easy to check that
${\bolds{\alpha}}^{\top}{\bolds{\Omega}}_{22.1}=0$, so
${\bolds{\alpha}}\perp\mbox{span}({\bolds{\Omega}}_{22.1})$ and
${\bolds{\alpha}}^{\top}{\bolds{\Omega}}_{22.1}^{+}{\bolds{\alpha}}=0$, and then
${\bolds{\alpha}}^{\top}({\mathbf{J}}^{\top}{\bolds{\Omega}}{\mathbf{J}})^{+}=0$. In this
situation,\break ${\mathbf{J}}({\mathbf{J}}^{\top}{\bolds{\Omega}}{\mathbf{J}})^{+}{\mathbf{J}}^{\top}\leq{\bolds{\Omega}}^{+}$ and the stick less-than sign holds
since ${\mathbf{J}}^{\top}{\bolds{\Omega}}{\mathbf{J}} \neq{\bolds{\Omega}}_{22.1}$ and
$1/{\bolds{\Omega}}_{11}>0$.
\hfill$\square$

\subsection{\texorpdfstring{Proof of Theorem \protect\ref{test}}{Proof of Theorem 2.2}}\label{sec4.5}

Under $H_0$, we can rewrite the index vector as
${\bolds{\beta}}=[\mathbf{e}\
{\mathbf{B}}]^{\top}(\sqrt{1-\|{\bolds{\omega}}^{(1)}\|^2},{\bolds{\omega}}^{(1)\tau})^{\top}$ where ${\mathbf{e}}=(1,0,\ldots,0)^{\top}$ is an $r$-dimensional vector,
\[
{\mathbf{B}}= \pmatrix{
{\bolds{0}}^{\top} & 0 \cr
{\mathbf{I}}_{r-1} & {\bolds{0}}
}
\]
is an $r\times(d-1)$ matrix and ${\bolds{\omega}}^{(1)}=(\beta_2,\ldots,\beta_r)^{\top}$ is an $(r-1)\times1$
vector. Let ${\bolds{\omega}}=(\sqrt{1-\|{\bolds{\omega}}^{(1)}\|^2},{\bolds{\omega}}^{(1)\top})^{\top}.$ So under $H_0$ the estimator is also the\vspace*{1pt}
local maximizer $\hat{\bolds{\omega}}$ of the problem
\[
\hat{Q}([\matrix{{\mathbf{e}}& {\mathbf{B}}}]^{\top}\hat{\bolds{\omega}})=\sup_{\|{\bolds{\omega}}^{(1)}\|< 1} \hat{Q}([\matrix{{\mathbf{e}}&
{\mathbf{B}}}]^{\top} {\bolds{\omega}}).
\]
Expanding $\hat{Q}({\mathbf{B}}^{\top}\hat{\bolds{\omega}})$ at $\hat{\bolds{\beta}}{}^{(1)}$
by a Taylor's expansion and noting that $\partial
\hat{Q}({\bolds{\beta}})/\break{\partial
{\bolds{\beta}}^{(1)}} |_{{\bolds{\beta}}^{(1)}=\hat{\bolds{\beta}}{}^{(1)}}=0$, then
$\hat{Q}(\hat{\bolds{\beta}})-\hat{Q}({\mathbf{B}}^{\top}\hat{\bolds{\omega}})=T_1+T_2+{\mo}_P(1)$, where
\begin{eqnarray*}
T_1 &=& -\frac{1}{2} \bigl(\hat{\bolds{\beta}}{}^{(1)}-{\mathbf{B}}^{\top}\hat{\bolds{\omega}} \bigr)^{\top} \frac{\partial^2 \hat{Q}({\bolds{\beta}})}{\partial
{\bolds{\beta}}^{(1)}\partial
{\bolds{\beta}}^{(1)\tau}} \bigg|_{{\bolds{\beta}}^{(1)}=\hat{\bolds{\beta}}{}^{(1)}}
\bigl(\hat{\bolds{\beta}}{}^{(1)}-{\mathbf{B}}^{\top}\hat{\bolds{\omega}} \bigr),\\
T_2 &=& \frac{1}{6} \bigl(\hat{\bolds{\beta}}{}^{(1)}-{\mathbf{B}}^{\top}\hat{\bolds{\omega}}
\bigr)^{\top}\\
&&{}\times\frac{\partial\{(\hat{\bolds{\beta}}{}^{(1)}-{\mathbf{B}}^{\top}\hat{\bolds{\omega}})^{\top} {\partial^2
\hat{Q}({\bolds{\beta}})}/{(\partial{\bolds{\beta}}^{(1)}\,\partial
{\bolds{\beta}}^{(1)\tau})}|_{{\bolds{\beta}}^{(1)}=\hat{\bolds{\beta}}{}^{(1)}}
(\hat{\bolds{\beta}}{}^{(1)}-{\mathbf{B}}^{\top}\hat{\bolds{\omega}}) \}
}{\partial
{\bolds{\beta}}^{(1)}}.
\end{eqnarray*}
Assuming the conditions in Theorem~\ref{asymptotic-normality} and
under the
null hypothesis $H_0$, it is easy to show that
\[
\sqrt{n}({\mathbf{B}}^{\top}\hat{\bolds{\omega}}-{\mathbf{B}}^{\top}{\bolds{\omega}
}) =
\frac{1}{\sqrt{n}}{\mathbf{B}}^{\top}{\mathbf{B}}({\mathbf{J}}^{\top} {\bolds{\Omega} \mathbf{J}})^{+} {\mathbf{G}}({\bolds{\beta}})+{\mo}_P(1).
\]
Combining this with (\ref{asym-expan-beta}), under the null hypothesis $H_0$,
%
\begin{eqnarray}\label{diff-expan}
&&\sqrt{n}\bigl(\hat{\bolds{\beta}}{}^{(1)}-{\mathbf{B}}^{\top}\hat{\bolds{\omega}}{}^{(1)}\bigr)\nonumber\\
&&\qquad =\frac{1}{\sqrt{n}}({\mathbf{J}}^{\top} {\bolds{\Omega} \mathbf{J}})^{1/2+}\{{\mathbf{I}}_{d-1}-({\mathbf{J}}^{\top} {\bolds{\Omega} \mathbf{J}})^{1/2}{\mathbf{B}}^{\top}{\mathbf{B}}({\mathbf{J}}^{\top}
{\bolds{\Omega} \mathbf{J}})^{1/2+}\}\\
&&\qquad \quad {}\times({\mathbf{J}}^{\top} {\bolds{\Omega} \mathbf{J}})^{1/2+}{\mathbf{G}}({\bolds{\beta}})+o_P(1).\nonumber
\end{eqnarray}
Since $\frac{1}{\sqrt{n}}{\mathbf{G}}({\bolds{\beta}})=\CO_P(1)$, $\frac
{\partial^2
\hat{Q}({\bolds{\beta}})}{\partial{\bolds{\beta}}^{(1)}\,\partial
{\bolds{\beta}}^{(1)\tau}} |_{{\bolds{\beta}}^{(1)}}=-n{\mathbf{J}}^{\top}{\bolds{\Omega} \mathbf{J}}
+{\mo}_P(n)$\vspace*{-3pt} and matrix ${\mathbf{J}}^{\top} {\bolds{\Omega} \mathbf{J}}$ has eigenvalues uniformly
bounded away from 0 and infinity,\vspace*{1pt} we have $\|\hat{\bolds{\beta}
}{}^{(1)}-{\mathbf{B}}^{\top}\hat{\bolds{\omega}}{}^{(1)}\|=\CO_P(n^{-1/2})$ and then
$|T_2|={\mo}_P(1)$.
Combining this and (\ref{diff-expan}), we have
\begin{eqnarray*}
\hat{Q}(\hat{\bolds{\beta}})-\hat{Q}({\mathbf{B}}^{\top}\hat{\bolds{\omega}
})&=&\frac{n}{2}
\bigl(\hat{\bolds{\beta}}{}^{(1)}-{\mathbf{B}}^{\top}\hat{\bolds{\omega}}{}^{(1)}\bigr)^{\top
}{\mathbf{J}}^{\top}{\bolds{\Omega} \mathbf{J}}\bigl(\hat{\bolds{\beta}}{}^{(1)}-{\mathbf{B}}^{\top}\hat
{\bolds{\omega}}^{(1)}\bigr)\\
&=& \frac{n}{2} {\mathbf{G}}^{\top}({\bolds{\beta}})({\mathbf{J}}^{\top}{\bolds{\Omega}
\mathbf{J}})^{1/2+}{\mathbf{P}}({\mathbf{J}}^{\top}{\bolds{\Omega} \mathbf{J}})^{1/2+}{\mathbf{G}}({\bolds{\beta}})
\end{eqnarray*}
with ${\mathbf{P}}={\mathbf{I}}_{d-1}-({\mathbf{J}}^{\top} {\bolds{\Omega} \mathbf{J}})^{1/2}{\mathbf{B}}^{\top}{\mathbf{B}}({\mathbf{J}}^{\top} {\bolds{\Omega} \mathbf{J}})^{1/2+}$. Here ${\mathbf{P}}$ is
idempotent having rank $d-r$, so it can be written as ${\mathbf{P}}= {\mathbf{S}}^{\top}{\mathbf{S}}$ where ${\mathbf{S}}$ ia a $(d-r)\times(d-1)$ matrix satisfying
${\mathbf{S}}{\mathbf{S}}^{\top}={\mathbf{I}}_{d-r}$. Consequently,
\begin{eqnarray*}
2\{\hat{Q}(\hat{\bolds{\beta}})-\hat{Q}({\mathbf{B}}^{\top}\hat{\bolds{\omega}})\}&=& \bigl(\sqrt{n} {\mathbf{S}}
({\mathbf{J}}^{\top}{\bolds{\Omega} \mathbf{J}})^{1/2+}
{\mathbf{G}}({\bolds{\beta}}) \bigr)^{\top} \bigl(\sqrt{n} {\mathbf{S}} ({\mathbf{J}}^{\top}{
\bolds{\Omega}
\mathbf{J}})^{1/2+} {\mathbf{G}}({\bolds{\beta}}) \bigr)\\
&\stackrel{\mathcal{L}} \longrightarrow&\chi^2(d-r).
\end{eqnarray*}
\end{appendix}\newpage

\section*{Acknowledgments}
The authors thank the Associate Editor and
two referees for their constructive comments and suggestions which
led to a great improvement over an early manuscript.

%

\begin{supplement}[id=suppA]
\stitle{Supplementary materials}
\slink[doi]{10.1214/10-AOS871SUPP} 
\sdatatype{.pdf}
\sfilename{aos871\_suppl.pdf}
\sdescription{Complete proofs of Proposition~\ref{prop-non}, (\ref
{eq-profile}) and (\ref{peqp}).}
\end{supplement}

%

\printaddresses

\end{document}